\documentclass[11pt]{article}

\usepackage{algorithm,xcolor}
\usepackage{algorithmic}

\usepackage{graphics}

\usepackage{subfigure}

\usepackage{graphicx}
\usepackage{rotating}

\title{Generalized Multiscale Finite Element Methods (GMsFEM)}
\author{Yalchin Efendiev\thanks{Department of Mathematics and ISC, Texas A \& M University, College Station, TX 77843} \and Juan Galvis\thanks{ISC, Texas A \& M University, College Station, TX 77843 and Departamento de Matem\'{a}ticas, Universidad Nacional de Colombia, 
Bogot\'a D.C., Colombia } \and Thomas Y. Hou\thanks{Computing and Mathematical Sciences, Pasadena, Caltech, CA 91125}}
\date{\today}

\usepackage{amsfonts}
\usepackage{amsmath,amssymb,amsthm}
\usepackage[english]{babel}

\newcommand{\RR}{\mathbb R}
\newtheorem{theorem}{Theorem}

\newtheorem{remark}[theorem]{Remark}

\usepackage{setspace}

\begin{document}

\selectlanguage{english}
\maketitle
%

\abstract{
In this paper, we propose a general approach called Generalized Multiscale
Finite Element Method (GMsFEM) for performing multiscale simulations
for problems without scale separation over a complex input space.
As in multiscale finite element methods (MsFEMs), the main idea
of the proposed approach is to construct a small dimensional local
solution space that can be used to generate an efficient and accurate 
approximation to the multiscale solution with a potentially high 
dimensional input parameter space. In the proposed approach,
we present a general procedure to construct the offline space 
that is used for a systematic enrichment of the coarse solution space
in the online stage. The enrichment in the online stage
is performed based on a spectral decomposition of the offline space.
In the online stage, for any input parameter, a multiscale 
space is constructed to solve the global problem on a coarse grid.
The online space is constructed via a spectral decomposition of
the offline space and by choosing the eigenvectors corresponding
to the largest eigenvalues. 
The computational saving is due to the fact that the construction
of the online multiscale space for any input parameter is fast and
this space can be re-used for 
solving the forward problem with
any forcing and boundary condition.
Compared with the other approaches where global snapshots are used, 
the local approach that we present in this paper allows us to eliminate 
unnecessary degrees of freedom on a coarse-grid level.
We present various examples in the paper and some numerical results
to demonstrate the effectiveness of our method.
}

\section{Introduction}

\subsection{Multiscale problems, the input-output relation,
 and the need for model reduction}

\begin{figure}[tb]
\centering
\includegraphics[width=4in, height=1.6in]{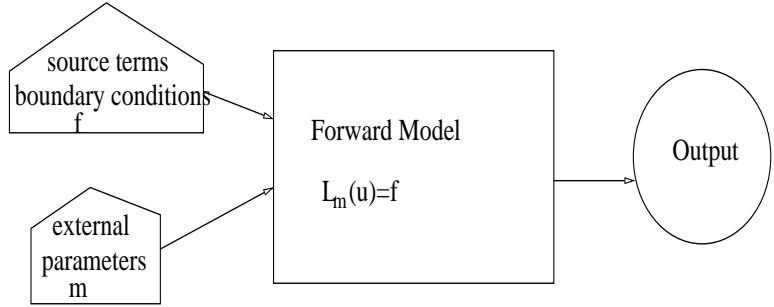}
\caption{Flow chart}
\label{scheme1}
\end{figure}

Many problems arising from various physical and engineering applications
 are multiscale in nature. Because of the presence of 
small scales and uncertainties 
in these problems, the direct simulations are prohibitively expensive.
Moreover, these problems are typically solved for many source terms 
with input parameter coming from a high dimensional parameter space.
For example, the flow in heterogeneous
porous media described by Darcy's equation is typically solved
for multiple source terms. Moreover, the permeability usually has
uncertainties which are parametrized in some sophisticated manner. 
In this case, one needs to solve many forward problems with different 
source terms and a wide range of permeabilities to make accurate 
predictions. These problems can be cast using an input-output relation 
(see Figure \ref{scheme1}) which is typically
done in reduced-order modeling. For the example of flow problems,
the input space consists of source terms and the permeability that
takes a value from  a large parameter space. 
The output space depends on the quantities
of interest and may consist of coarse-grid solutions or some other 
integrated quantities with respect to the solution. In many applications,
the output space is typically smaller than the input space. The design of  a
general multiscale finite element framework that takes advantage of 
the effective low dimensional solution space for multiscale problems 
with high dimensional input space is the main objective of this paper.

Due to a large number of forward simulations, the computational
effort can be tremendous to learn and process the output space 
given the high dimensionality of the input parameter space. In many 
of these problems, the solution space can be approximated by a low 
dimensional manifold via some model reduction tools. The main objective
of reduced-order models is to represent the solution space with a 
small dimensional space. However, many existing reduced-order methods fail to give a small dimensional output solution space when the physical 
solution has multiscale structures. Another major limitation of the 
current reduced-order methods is that the reduced solution space needs 
to be regenerated for different forcing or boundary conditions. 
The general multiscale finite element framework proposed in this paper
is designed to remove these two limitations by dividing the construction
of reduced basis into the offline and online steps, and constructing
our online multiscale bases from a reduced localized offline 
solution space.

\subsection{Local and global model reduction concepts}

Many local, global, and
local-global model reduction techniques have been developed. 
The main idea
of these methods is to find a small dimensional space that can represent
the solution space given the input space. 

Global model reduction techniques 
(see e.g., 
\cite{LSBook_Joao, LSBook_Chen, AntoulasBook, volkwein05})
construct a space of global
fields that can approximate the solution space. 
One can, for example, consider a space of exhaustive global snapshots
obtained by solving the global problem for many input parameters.
This space can be further reduced using a spectral decomposition.
In practice, the resulting space is constructed
by solving global problems for some selected input parameters, right
hand sides, and boundary conditions. These
methods have been used with some success in practice. However, when 
the right hand sides or boundary conditions are changed, the resulting
reduced space must be recomputed.

\begin{figure}[tbp]
\centering
\includegraphics[width=5.5in, height=3.in]{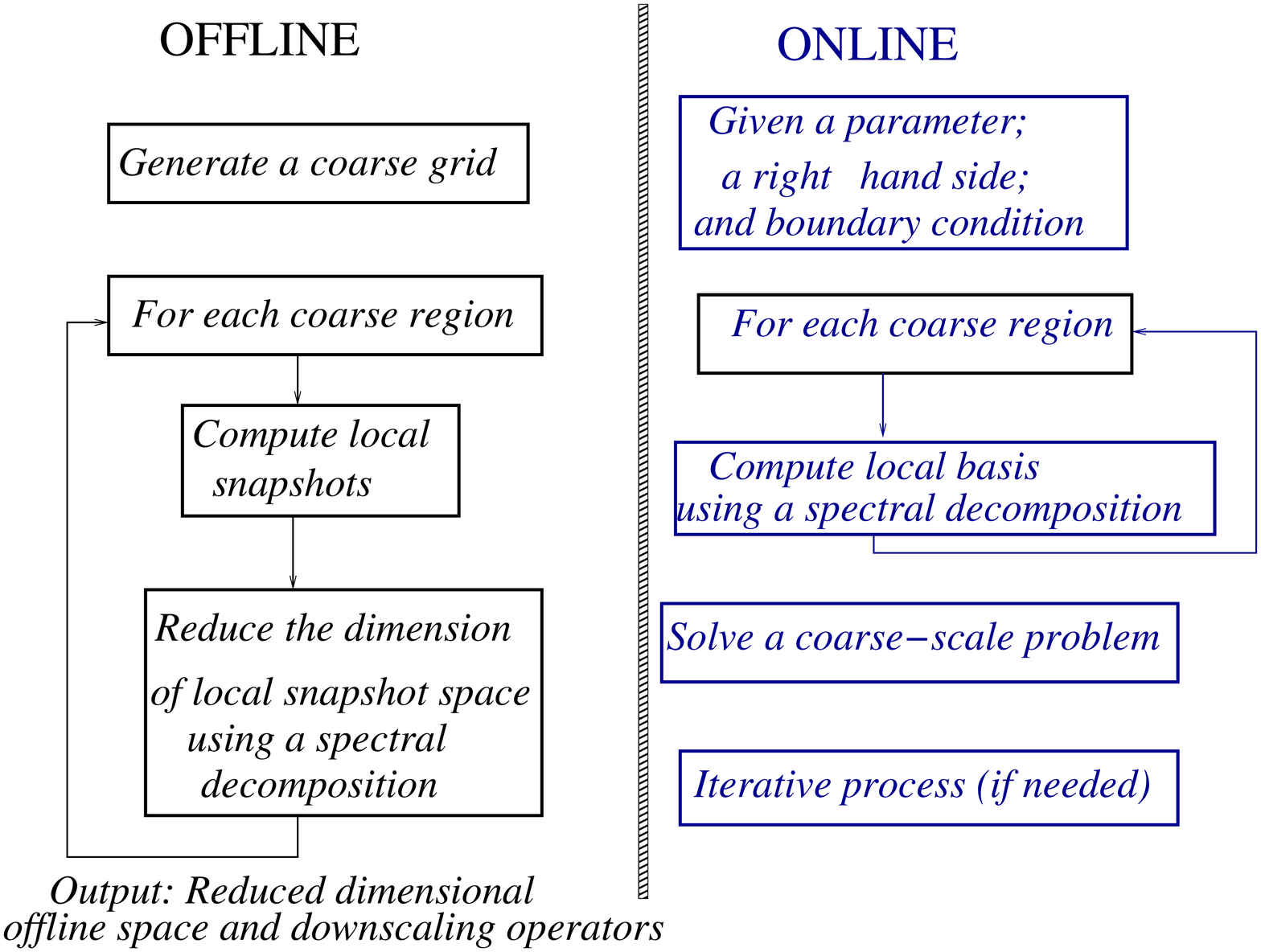}
\caption{Flow chart}
\label{scheme}
\end{figure}

Local approaches 
(e.g., see \cite{hw97, krog11, aarnes04, aej07, AKL, arbogast02, apwy07, cdgw03, ch03, eghe05, jennylt03} for upscaling and multiscale methods)
attempt to approximate 
the solution in local (coarse-grid) regions for all input parameters
without computing global snapshots of solutions. 
Local approaches first compute an offline space (possibly small 
dimensional)
which is used to compute multiscale basis functions at the online stage.
The local approximation space at the online stage is computed by
finding a subspace of offline space for a given input parameter
(see Figure \ref{scheme}).

Local approaches
can be effective as they avoid the computation of global snapshots.
Local approaches become more effective if the restriction of
the solution space onto a local region has a small dimension.
This is the case if the dimension of the space of solutions restricted to
a coarse region is smaller than the dimension of the fine-grid
space 
within this coarse region.
For example, if the parameter is a coarse-grid scalar function, then
at the coarse-grid level, this parameter is a scalar. While if we
consider this problem from the point of view of a global model 
reduction, then the parameter belongs to a large dimensional space 
and this may not be amenable to computations.

One of advantages of local approaches is that they eliminate the
unnecessary degrees of freedom in the parameter space
 at the coarse-grid level. In global methods,
one first needs to compute many expensive global snapshots 
and many snapshots may not contribute to the solution at the 
online stage. In local approaches, these values of the parameter 
are identified at the coarse level inexpensively. Moreover, local 
approaches can easily handle large-scale parameter
space when the parameter is a coarse-grid function and local
approximation spaces are usually independent of the source
terms or boundary conditions. We will further elaborate these 
issues in the paper.

\subsection{This paper}

In this paper, 
we introduce a general multiscale framework, which we call the 
Generalized Multiscale
Finite Element Method (GMsFEM). This method
incorporates complex input space information
and the input-output relation. It systematically enriches the 
coarse space through our local construction. Our approach, as in many
multiscale and model reduction techniques, divides the computation into
two stages: offline; and online. In the offline stage, we
construct a small dimensional space that can be efficiently
used in the online stage to construct multiscale basis functions.
These multiscale basis functions can be re-used for any input parameter
to solve the problem on a coarse-grid. Thus, this provides a substantial
computational saving at the online stage. Below, we present an outline 
of the algorithm
and a chart that depicts our algorithm in Figure \ref{scheme}. 

\begin{itemize}
\item[1.]  Offline computation:
\begin{itemize}
\item 1.0. Coarse grid generation;  
\item 1.1. Construction of snapshot space that will be used to compute an offline space. 
\item 1.2. Construction of a small dimensional offline space by performing dimension reduction in the space of global snapshots.
\end{itemize}
\item[2.] Online computations:
\begin{itemize}
\item 2.1. For each input parameter, compute multiscale basis functions;
\item 2.2. Solution of a coarse-grid problem for any force term and boundary condition;
\item 2.3. Iterative solvers, if needed.
\end{itemize}
\end{itemize}

In the offline computation,
we first set up a coarse grid where each coarse-grid block consists
of a connected union of fine-grid
blocks. The construction of snapshot space in Step 1.1 involves solving local
problems for various choices of input parameters. 
This space is used to construct the offline space in Step 1.2 via a
spectral decomposition of the snapshot space. The 
snapshot space in a coarse region
 can be replaced by
the fine-grid space associated with this coarse space; however, in many
applications, one can judiciously choose the space of snapshots to avoid 
expensive offline space construction. The offline space in Step 1.2.
is constructed by spectrally decomposing the space of snapshots.
This spectral decomposition is typically based on the offline eigenvalue
problem. The spectral decomposition enables us to select the 
high-energy elements from the offline space by choosing those eigenvectors 
corresponding to the largest eigenvalues.
More precisely, we seek a subspace of the snapshot space such
that it can approximate any element of the snapshot space in 
the appropriate sense defined via auxiliary bilinear forms. 

In the online step 2.1, for a given input parameter, we compute
the required online coarse space. In general, we want this to be
a small dimensional subspace of the offline space. This space is
computed  by performing a spectral decomposition in the offline space
via an eigenvalue problem. Furthermore, the eigenvectors corresponding
to the largest eigenvalues are identified and used to form the online 
coarse space. The online coarse space is used within the finite element
framework to solve the original global problem. Here, we propose several
options such as the Galerkin coupling of multiscale basis functions,
the Petrov-Galerkin coupling of multiscale basis functions, etc.
 In some of these coupling approaches, the choice of the initial partition 
of unity (that can be computed in the offline or online stage)
 is important and it will be discussed in the paper.

Our techniques differ from many previous approaches that are based
on the homogenization theory. In the homogenization based methods,
one usually constructs local approximation based on  local
solves and these approaches do not provide a systematic procedure
to complement the local spaces. It is important to note that one needs
to systematically complement the local spaces in order to converge to
the fine-grid solution. How to develop an online systematic enrichment
procedure and how to construct the 
initial partition of unity functions
 play a crucial role in obtaining a low dimensional offline
space. These issues are central points of our proposed method. 

We also discuss iterative solvers that use the coarse spaces and 
iterate on the residual to converge to the fine-scale solution. 
These iterative solvers serve as an online correction of the
coarse-grid solution. We consider two-level domain decomposition 
preconditioners and some other methods where the importance of 
appropriately chosen multiscale coarse spaces has been demonstrated
in the literature \cite{tw, Tarekbook}.
We will discuss how the choice of
coarse spaces  yields optimal iterative solvers where the number
of iterations is independent of the high contrast in the media
properties.

\section{A Generalized Multiscale Finite Element Method.}
\label{prelim}

To describe the GMsFEM for linear problems, we consider
\begin{equation}
\label{main:red}
L_\mu(u)=f,
\end{equation}
subject to some boundary conditions, where $\mu$ is the parameter. 
For example,
\begin{equation}
\label{eq:ellipcase}
L_\mu(u)=-\mbox{div}(\kappa(x;\mu)\nabla u).
\end{equation}
Here, the operator $L$ may depend on various spatial
fields, e.g., heterogeneous conductivity fields, convection fields, 
reaction fields, and so on.
The dependence of the solution from these fields is nonlinear, 
while the solution linearly depends
on external source terms $f$ and boundary conditions.
We assume there is a bilinear
form associated with  the operator $L$ that allows us to write 
the variational form of  (\ref{main:red}) in the form
\begin{equation}
\label{eq:bilinearD}
{\kappa}(u,v;\mu)=\langle L_\mu(u), v\rangle,
\end{equation}
for all test function $v$ and 
where ${\kappa}(\cdot,\cdot;\mu)$ is bilinear, coercive, and continuous
for each $\mu$, and $\langle\cdot,\cdot\rangle$
is an inner product.
 We assume that ${\kappa}(\cdot,\cdot;\mu)$
is sufficiently smooth with respect to $\mu$.

The proposed algorithms have advantages when $L_\mu(u)$ has an affine
representation defined as follows:
\begin{equation}\label{hyp}
L_\mu(u) := \sum\limits_{q=1}^Q \Theta_q(\mu)L_q(x).
\end{equation}
Here, $L_q(u)$ are heterogeneous operators 
with multiple scales with the coefficients that have high contrast,
the parameter $\mu \in \Lambda \subset \RR^p$ is possibly a coarse-grid 
function and the functions $\Theta_q : \Lambda \to \RR$. 
In terms of the corresponding bilinear form, we have
\[
{\kappa}(u,v;\mu)=\sum\limits_{q=1}^Q \Theta_q(\mu) \kappa_q(u,v).
\]
This affine representation allows pre-computing coarse-scale projections
of $L_q(x)$ in the offline stage and 
using them in the online stage. This reduces 
the computational cost considerably.

Before discussing the GMsFEM, we introduce the notion of a coarse grid.
Let $\mathcal{T}^H$  be a usual  conforming partition of $D$ into finite elements
 (triangles, quadrilaterals, tetrahedrahals,...). We call this partition 
the coarse grid and 
assume that this coarse grid is partitioned into fine-grid blocks. 
Each coarse-grid block is a connected union of fine-grid blocks.
We denote by $N_v$ the number of coarse nodes, by $\{x_i\}_{i=1}^{N_v}$ the vertices of 
the coarse mesh $\mathcal{T}^H$ and define the neighborhood of the node $x_i$ by
\begin{equation}\label{eq:def:omegai}
\omega_i=\bigcup\{ K_j\in\mathcal{T}^H; ~~~ x_i\in \overline{K}_j\}
\end{equation}
(see Figure \ref{coarse_fine_ovs}).

The GMsFEM has a structure similar to that of MsFEM. The main difference
between the two approaches is that we systematically enrich coarse 
spaces in GMsFEM and generalize it by considering an 
input space consisting of parameters and source terms.
 In the first step of GMsFEM (offline stage), we
 construct the space of `snapshots', $V_{\text{snapshots}}^{\omega_i}$, 
a large dimensional space of local solutions. In the next step of the 
offline computation, we reduce the space $V_{\text{snapshots}}^{\omega_i}$ 
via some spectral procedure to $V_{\text{off}}^{\omega_i}$.
In the second stage (online stage), for each input parameter,
we construct a corresponding local space, $V_{\text{on}}^{\omega_i}$ that is used
to solve the problem at the online stage for the given input parameter.
Our systematic
approach allows us to increase the dimension of the coarse space 
and achieve a convergence.



\subsection{An illustrative example}
\label{sec:ill_example}

Before presenting details of GMsFEM, we will present a simple
example demonstrating the 
main concept of GMsFEM. We consider
\[
-\mbox{div}(\kappa(x;u)\nabla u)=f\ \text{in}\ D,
\]
$u=0$ on $\partial D$ and assume $u_0\leq u(x) \leq u_N$,
where $u_0$ and $u_N$ are pre-defined constants.
We assume the interval $[u_0,u_N]$ is divided
into $N$ equal regions $u_0<u_1<...<u_{N-1}<u_N$.

As for the space of {\it snapshots}, $V_{\text{snapshots}}^{\omega_i}$,
we can consider
\[
-\mbox{div}(\kappa(x;u_j)\nabla \psi_{l,j}^{\text{snap}})=0\ \ \text{in}\ \omega_i,
\]
$ \psi_{l,j}^{\text{snap}}=\delta_l(x)$ on $\partial \omega_i$. Here,
$\delta_l(x)$ are some set of functions defined on $\partial \omega_i$,
e.g., unit source terms.  We
can also consider the space of fine-grid functions within $\omega_i$ 
as the space of  {snapshots}.

As for {\it offline} space, $V_{\text{off}}^{\omega_i}$,
 we perform a spectral decomposition
of the space of snapshots. We re-numerate
the snapshot functions in $\omega_i$ by $\psi_l^{\text{snap}}$. We consider
\[
(S^{\text{off}}:=)s_{mn}=\int_{\omega_i} \sum_l t_l \kappa(x;u_l) \nabla \psi_m^{\text{snap}}\cdot \nabla \psi_n^{\text{snap}},\ \ 
(A^{\text{off}}:=)a_{mn}=\int_{\omega_i} \sum_lt_l  \kappa(x;u_l) \psi_{m}^{\text{snap}} \psi_{n}^{\text{snap}},
\]
where $t_l$ are some weights.
Here, $S^{\text{off}}=(s_{mn})$, and $A^{\text{off}}:=(a_{mn})$.
The choices for the matrix $A^{\text{off}}$ will be discussed later.
To generate the offline space, $V_{\text{off}}^{\omega_i}$, we choose
the largest $M_{\text{off}}$ eigenvalues (see later discussions on
 $M_{\text{off}}$)
of
\[
A^{\text{off}} \Psi_m ^{\text{off}}= \lambda_m^{\text{off}} S^{\text{off}} \Psi_m^{\text{off}}
\]
and find the corresponding eigenvectors in the space of $V_{\text{snapshots}}^{\omega_i}$
by multiplication,  $\sum_{j}\Psi_{ij}^{\text{off}} \psi^{\text{snap}}_j$, where
$\Psi_{ij}^{\text{off}}$ are coordinates of the vector $\Psi_i^{\text{off}}$.
Note that the offline space is computed across all $u_0$,..., $u_N$. 
More precisely, if we reorder the snapshot 
functions using a single index based on the decay of eigenvalues to create the matrices
\[
R_{\text{snap}} = \left[ \psi_{1}^{\text{snap}}, \ldots, \psi_{M_{\text{snap}}}^{\text{snap}} \right],
\]
then
\begin{equation*}
 \displaystyle S^{\text{off}}= [s^{\text{off}}_{mn}] = \int_{\omega_i} \overline{\kappa}(x) \nabla \psi_m^{\text{snap}} \cdot \nabla \psi_n^{\text{snap}} = R_{\text{snap}}^T \overline{S} R_{\text{snap}},
 \end{equation*}
\begin{equation*}
 \displaystyle A^{\text{off}} = [a^{\text{off}}_{mn}] = \int_{\omega_i} {\overline{\kappa}}(x) \psi_m^{\text{snap}} \psi_n^{\text{snap}} = R_{\text{snap}}^T \overline{A} R_{\text{snap}},
 \end{equation*}
where $\overline{S}$ and $\overline{A}$ denote  fine-scale matrices
corresponding to the stiffness and mass matrices with the permeability
$\overline{\kappa}(x) =  \sum_l t_l \kappa(x;u_l)$.
One can use modified $\overline{\kappa}$ in the computation
of the mass matrix (see \cite{egt11}).

To compute the {\it online} space for a given $u_q$
(the value around which the global problem is linearized),
 we consider
an eigenvalue problem in $V_{\text{off}}^{\omega_i}$. Let us denote
the basis of  $V_{\text{off}}^{\omega_i}$ by $\psi_m^{\text{off}}$. We consider
a spectral decomposition of $V_{\text{off}}^{\omega_i}$ via
\[
S^{\text{on}}=(s_{mn})=\int_{\omega_i}  \kappa(x;u_q) \nabla \psi_m^{\text{off}}\cdot \nabla \psi_n^{\text{off}},\ \ 
A^{\text{on}}=(a_{mn})=\int_{\omega_i}  \kappa(x;u_q) \psi_{m}^{\text{off}} \psi_{n}^{\text{off}}.
\]
To generate the online space, $V_{\text{on}}^{\omega_i}$, we choose
the largest $M_{\text{on}}$ eigenvalues of
\[
A^{\text{on}} \Psi_m^{\text{on}} = \lambda_m^{\text{on}}S^{\text{on}} \Psi_m^{\text{on}}
\]
and find the corresponding eigenvectors 
in the space of $V_{\text{off}}^{\omega_i}$
by multiplication,  $\sum_{j}\Psi_{ij}^{\text{on}} \psi^{\text{off}}_j$, where
$\Psi_{ij}^{\text{on}}$ are coordinates of the vector $\Psi_i^{\text{on}}$,
denote these basis functions by $\psi_m^{\text{on}}$. 
More precisely, if we reorder the offline
functions using a single index based on the decay of eigenvalues to create the matrices
\[
R_{\text{off}} = \left[ \psi_{1}^{\text{off}}, \ldots, \psi_{M_{\text{off}}}^{\text{off}} \right]
\]
that are defined on the fine grid,
then
\begin{equation*}
 \displaystyle S^{\text{on}}= [s^{\text{on}}_{mn}] = \int_{\omega_i} \kappa(x; u_q) \nabla \psi_m^{\text{off}} \cdot \nabla \psi_n^{\text{off}} = R_{\text{off}}^T {S} R_{\text{off}},
 \end{equation*}
\begin{equation*}
 \displaystyle A^{\text{on}} = [a^{\text{on}}_{mn}] = \int_{\omega_i} {\kappa}(x; u_q) \psi_m^{\text{off}} \psi_n^{\text{off}} = R_{\text{off}}^T {A} R_{\text{off}},
 \end{equation*}
where ${S}$ and ${A}$ denote  fine scale matrices
corresponding to the stiffness and mass matrices.

At the final stage, these basis functions $ \Psi_m^{\text{on}}$
in each $\omega_i$ will be coupled via a global
formulation, e.g., Galerkin formulation. 
In this case, the eigenfunctions are multiplied by the partition
of unity functions to obtain a conforming basis.
In this nonlinear example,
one can use an iterative Picard iteration at the previous
value of the solution $u^n(x)$, and in each iteration, a global
problem is solved with $V_{\text{on}}$ for the value of $u^n(x)$
averaged over a coarse block.

\subsection{Step 1. Local multiscale basis functions (offline stage)}

{\bf 1.1. Generating snapshots}

We call the local spatial fields (or local solutions ) 
 as ``local snapshots''. These 
local spatial fields are used to construct the offline space. 
This space consists of 
spatial fields defined on a fine grid, i.e., they are
vectors of the dimension of fine-grid resolution of the coarse region.
A common option for local snapshots
is to use a fine-grid space (e.g., fine-grid linear functions)
 that resolves the coarse-grid block.
However, in some cases, one can consider smaller and more appropriate
spaces to construct the local snapshot space. 
This will be discussed next.

In the first step, a
local reduced-order approximation is constructed based on the input space.
Here, we will discuss two approaches and emphasize only the first 
approach where we will construct local approximate solutions by
 assuming that source term $f$ is a smooth function. In this 
construction, the source term will not enter in the space
of local solutions and its effect will be captured at the coarse-grid 
level via a global coupling. 

\begin{remark}
One can often ignore lower order
terms in the operator $L$  and use an operator
different from the original one
on a coarse grid to construct snapshots. To formalize this step, we
assume that there exists $\widetilde{L}_\mu$ 
such
that if
\[
\widetilde{L}_\mu(\widetilde{u})=0 \ \ \text{and}\ \widetilde{u}=u\ \text{\text{on}}\ \partial K
\]
then
\[
\|\widetilde{u}-u\|\preceq\ \delta_H,
\] 
where $\delta_H\rightarrow 0$ as $H\rightarrow 0$.
In this case, we have a corresponding bilinear form $\widetilde{a}(u,v;\mu)$.
For example, if 
$L(u)=-\mbox{div}(\kappa(x/\epsilon;\mu)\nabla u) + q(x;\mu)\cdot \nabla u$, one can use
$ \widetilde{L}(u)=-\mbox{div}(\kappa(x/\epsilon;\mu)\nabla u) $ (see e.g., 
\cite{blp78}), provided, e.g., $q$ is a bounded function.
To avoid a cumbersome notation, we do not use $\widetilde{L}$ in the rest 
of the presentation.


\end{remark}

\begin{figure}[tbp]
\centering
\includegraphics[width=3in, height=2in]{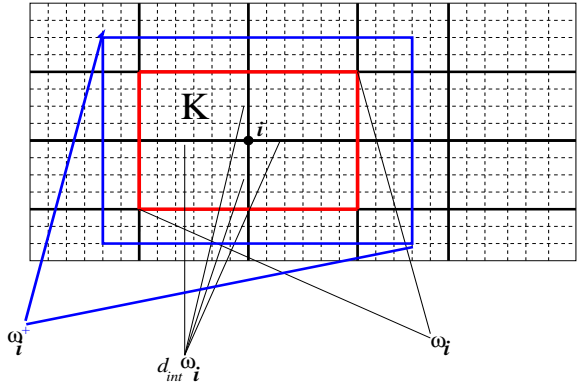}
\caption{Schematic description of oversampled regions.} 
\label{coarse_fine_ovs}
\end{figure}

We construct local snapshots of the solutions that approximate the space
of solutions generated by
\begin{equation}
\label{eq:tilde}
{L}_\mu(v)=0
\end{equation}
in each subdomain $\omega_i$ (see Figure \ref{coarse_fine_ovs}).
We will consider two choices, though one can consider other options. 

{\bf First choice.} We consider snapshots generated by
\begin{equation}
\label{eq:case1.1}
{L}_{\mu_j}(\psi_{l,j}^{\omega_i})=0\ \text{in} \ \omega_i
\end{equation}
with boundary conditions
\begin{equation}
\label{eq:case1.2}
\psi_{l,j}^{\omega_i}=b_l\ \text{in} \ \partial \omega_i,
\end{equation}
with $b_l$ 
being selected shape or  basis functions  along the boundary $\partial\omega_i$.

One can also use Neumann boundary conditions or
boundary conditions $b_l$ defined on larger domains for generating snapshots.

{\bf Second choice.} 
We can use local fine-scale spaces consisting of fine-grid basis
functions within a coarse region. In this case, the offline spaces
will be computed on a fine grid directly.
The local fine-grid space has an advantage  if the dimension of the local 
fine spaces  is comparable to the dimension of  
$V_{\text{snapshots}}^{\omega_i}$ computed by solving local problems
as in the first choice.
We use local fine-grid basis functions as a snapshot space in \cite{egt11}
and our earlier works \cite{ge09_2} which can be difficult to extend, in general.


One can also construct local snapshots by solving the local spectral
problem 
\[
\widetilde{A}_{\mu_j}(\psi_{l,j}^{\omega_i})=\lambda_l \widetilde{S}_{\mu_j}(\psi_{l,j}^{\omega_i})\ \text{in} \ \omega_i
\]
with homogeneous Neumann boundary conditions and
for some operators $\widetilde{A}_{\mu_j}$ and $\widetilde{S}_{\mu_j}$
 and selecting the dominant 
eigenvectors.



\begin{remark}
{\bf On generating snapshots via the approximation of source term $f$.} 
The snapshot spaces generated above can be
larger than the space corresponding 
to local snapshots obtained from 
\[
L_\mu(v_f)=f
\]
where $f$ runs over the input space corresponding to the source term.
Here, we take the restriction of $v_f$ in $\omega_i$ to generate
the space of local snapshots. 
This space can be taken as a span of $\phi_i^u$,
where $\phi_i^u$ are restrictions of the solutions of 
\[
L_\mu(\phi_i^u)=\phi_i^f
\]
onto a coarse region $\omega_i$.
Here $\phi_i^f$ are local basis functions that represent
\[
f=\sum_i f_i \phi_i^f.
\]
These basis functions are global; however, they can be localized at  a 
coarse-grid level (see e.g., Owhadi and Zhang, 2011, \cite{oz11}) and 
we need a further dimension 
reduction of the space following Step 1.2. One can use the restriction
of these global snapshots on the boundary of the coarse-grid block
to compute the space of snapshots. This space of snapshots is
an appropriate space of functions for which spectral decomposition
needs to be performed.
\end{remark}


{\bf 1.2. Generating the offline space}

Once the space of snapshots, $V_{\text{snapshots}}^{\omega_i}=Span(\psi^{\omega_i}_{l,j})$, is constructed
for each $\omega_i$, spectral approaches are needed to orthogonalize and
possibly
reduce the dimension of this space. As a result, we will obtain 
the offline space, $V_{\text{off}}^{\omega_i}$ that will be used to construct
multiscale basis functions in the online stage.

To perform a dimension reduction, we consider an auxiliary
spectral decomposition of the space 
$V_{\text{snapshots}}^{\omega_i}=Span(\psi^{\omega_i}_{l,j})$
for each $\omega_i$. Our objective is to construct a 
possibly small
dimensional space  $V_{\text{off}}^{\omega_i}$ and use it for
constructing multiscale basis functions for each $\mu$
in the online stage. In general, we will seek
the subspace $V_{\text{off}}^{\omega_i}$ such that for any
$\mu$ and
$\psi\in V_{\text{snapshots}}^{\omega_i}(\mu)$ (
$V_{\text{snapshots}}^{\omega_i}(\mu)$ is the space of
snapshots which are computed for a given $\mu$), there exists 
$\psi_0\in V_{\text{off}}^{\omega_i}$, such that, for all $\mu$,
\begin{equation}
\label{eq:off1}
a_{\omega_i}^{\text{off}}(\psi-\psi_0,\psi-\psi_0;\mu)\preceq {\delta} 
s_{\omega_i}^{\text{off}}(\psi-\psi_0,\psi-\psi_0;\mu),
\end{equation}
where $a_{\omega_i}^{\text{off}}(\phi,\phi;\mu)$ and 
$s_{\omega_i}^{\text{off}}(\phi,\phi;\mu)$
are auxiliary bilinear forms. 
In computations, this involves solving an eigenvalue problem
with a mass matrix  and the basis functions are selected based on
dominant eigenvalues. Note that this eigenvalue problem is formed
in the snapshot space. 
We will discuss two procedures for constructing
$V_{\text{off}}^{\omega_i}$.


\begin{remark}
\label{rem:a_vs_s}
In general, $a_{\omega_i}^{\text{off}}$ and $s_{\omega_i}^{\text{off}}$ contain
partition of unity functions, penalty terms, and other discretization
factors that appear in coarse-grid finite element formulations.
The norm corresponding to  $s_{\omega_i}^{\text{off}}$ needs to be stronger,
in general, to allow the decay of eigenvalues.

\end{remark}


We consider two other options.

{\bf Option 1.} We consider 
 $a_{\omega_i}^{\text{off}}(\phi,\phi;\mu)$ and 
$s_{\omega_i}^{\text{off}}(\phi,\phi;\mu)$ to be independent of 
$\mu$, i.e.,
\[ 
a_{\omega_i}^{\text{off}}(\phi,\phi;\mu)=a_{\omega_i}^{\text{off}}(\phi,\phi)\ \ \ 
\text{and}\ \ \ 
s_{\omega_i}^{\text{off}}(\phi,\phi;\mu)=s_{\omega_i}^{\text{off}}(\phi,\phi).
\]
 In this case, finding $V_{\text{off}}^{\omega_i}$ reduces to 
performing spectral decomposition of $ V_{\text{snapshots}}^{\omega_i}$
with corresponding inner products.
A subspace $V_{\text{off}}^{\omega_i}$ of  
$V_{\text{snapshots}}^{\omega_i}$ such that for each
$\psi\in V_{\text{snapshots}}^{\omega_i}$, 
there exists $\psi_0\in V_{\text{off}}^{\omega_i}$ such that
\begin{equation}
\label{eq:eigineq}
a_{\omega_i}^{\text{off}}(\psi-\psi_0,\psi-\psi_0)\preceq \delta s_{\omega_i}^{\text{off}}(\psi-\psi_0,\psi-\psi_0)
\end{equation}
 for some prescribed error tolerance $\delta$.
We give some examples for the elliptic equation 
(\ref{eq:ellipcase}).

\noindent {\bf Example.}
In this example, we 
average
the parameter $\mu$ to obtain
$a^{\text{off}}$ and $s^{\text{off}}$.
\begin{equation}
\begin{split}
s_{\omega_i}^{\text{off}}(\psi,\psi)=\sum_{j} t_j a_{\omega_i}(\psi,\psi;\mu_j)=
\sum_j t_j \int_{\omega_i} \kappa(x;\mu_j) |\nabla \psi|^2\\ 
a^{\text{off}}_{\omega_i}(\psi,\psi)=\sum_{j}t_j\int_{\omega_i}
\kappa(x,\mu_j)|\nabla(\chi_i\psi)|^2, \ \ \text{or} \ \ 
a^{\text{off}}_{\omega_i}(\psi,\psi)=\sum_{j}t_j\int_{\omega_i}
\kappa(x,\mu_j)|\psi|^2,
\end{split}
\end{equation}
where $t_j$ are non-negative weights (a case with a fixed value of
$\mu$ is a special case) and $\chi_i$ is a partition of unity function
supported in $\omega_i$.

To formulate the eigenvalue problem corresponding to (\ref{eq:eigineq}),
we will re-numerate the basis,
$V_{\text{snapshots}}^{\omega_i}=Span(\psi_{J}^{\omega_i})$. Then
(\ref{eq:eigineq}) will yield an algebraic eigenvalue problem
\[
A^{\text{off}}\Psi_l^{\text{off}}=\lambda_l^{\text{off}}S^{\text{off}} \Psi_l^{\text{off}}.
\]
Here,
$A^{\text{off}}=(A^{\text{off}}_{IJ})$ and $S^{\text{off}}=(S^{\text{off}}_{IJ})$ are corresponding
local matrices that are given by 
\[
A^{\text{off}}_{IJ}=a_{\omega_i}^{\text{off}}(\psi_I^{\omega_i},\psi_J^{\omega_i}),
\quad 
S^{\text{off}}_{IJ}=s_{\omega_i}^{\text{off}}(\psi_I^{\omega_i},\psi_J^{\omega_i}),
\]
$\psi_I^{\omega_i}, \psi_J^{\omega_i}\in V_{\text{snapshots}}^{\omega_i}$. Here, we assume 
that $S^{\text{off}}_{ij}$  is a non-degenerate positive definite matrix.
Note that the matrices $A^{\text{off}}$ and $S^{\text{off}}$ are
computed in the snapshot space.
The space $V_{\text{off}}^{\omega_i}$ is constructed by selecting 
$L_i$ eigenvectors corresponding to largest
eigenvalues. If we assume that 
\[
\lambda_1^{\text{off}}\geq ...\geq \lambda_N^{\text{off}},
\]
then we choose the first $L_i$ eigenvectors
in each $\omega_i$ corresponding to the largest 
$L_i$ eigenvalues to form the offline space. 
In particular,
the corresponding eigenvectors 
are computed 
by multiplication,  $\sum_{j}\Psi_{ij}^{\text{off}} \psi^{\omega_l}_j$
in each $\omega_l$,  where
$\Psi_{ij}^{\text{off}}$ are coordinates of the vector $\Psi_i^{\text{off}}$.

More precisely, if we reorder the snapshot
functions using a single index based on the decay of eigenvalues to create the matrices
\[
R_{\text{snap}} = \left[ \psi_{1}^{\text{snap}}, \ldots, \psi_{M_{\text{snap}}}^{\text{snap}} \right],
\]
then
\begin{equation*}
 \displaystyle S^{\text{off}}= [s^{\text{off}}_{mn}] = \int_{\omega_i} \overline{\kappa}(x) \nabla \psi_m^{\text{snap}} \cdot \nabla \psi_n^{\text{snap}} = R_{\text{snap}}^T \overline{S} R_{\text{snap}},
 \end{equation*}
\begin{equation*}
 \displaystyle A^{\text{off}} = [a^{\text{off}}_{mn}] = \int_{\omega_i} {\overline{\kappa}}(x) \psi_m^{\text{snap}} \psi_n^{\text{snap}} = R_{\text{snap}}^T \overline{A} R_{\text{snap}},
 \end{equation*}
where $\overline{S}$ and $\overline{A}$ denote  fine-scale matrices
corresponding to the stiffness and mass matrices with the permeability
$\overline{\kappa}(x)=\sum_j t_j \kappa(x;\mu_j)$, which is selected independent of $\mu$.

{\bf Option 2.}
The second option for a numerical setup of  
(\ref{eq:off1})
 is the following (see \cite{egt11}). First, we will compute the local
spaces for each $\mu$ (from a pre-selected exhaustive set), 
$V_{\text{off},\mu}^{\omega_i}(\mu)$ such that
for any $\mu$ in this set and
$\psi\in V_{\text{snapshots}}^{\omega_i}(\mu)$ (i.e.,
snapshots are computed for a given $\mu$), there exists 
$\psi_0\in V_{\text{off},\mu}^{\omega_i}(\mu)$, such that
\begin{equation}
\label{eq:off2}
a_{\omega_i}^{\text{off}}(\psi-\psi_0,\psi-\psi_0;\mu)\preceq {\delta} 
s_{\omega_i}^{\text{off}}(\psi-\psi_0,\psi-\psi_0;\mu).
\end{equation}
This involves finding the dominant eigenvectors for selected
set of $\mu$'s.
Secondly, we seek a small dimensional space 
$V_{\text{off}}^{\omega_i}$ (across all $\mu$'s) such that any
$\psi_0\in V_{\text{off},\mu}^{\omega_i}(\mu)$ can be approximated by an 
element of $V_{\text{off}}^{\omega_i}$. One can do it by defining a distance
function between two spaces $ V_{off\mu}^{\omega_i}(\mu_1)$
and  $ V_{\text{off},\mu}^{\omega_i}(\mu_2)$ and identifying a small
number of $\mu$'s such that the spaces  $V_{\text{off},\mu}^{\omega_i}(\mu)$
corresponding to these $\mu$'s approximate 
the space spanned by  $V_{\text{off},\mu}^{\omega_i}(\mu)$ for all $\mu$'s
(see \cite{egt11}).
Then, the elements of these  $V_{\text{off},\mu}^{\omega_i}(\mu)$'s 
will form the space $V_{\text{off}}^{\omega_i}$.

Note that
it is important to have a small dimensional space of snapshots
to reduce the computational cost that arises in
the calculations of multiscale basis functions.
The result of Step 2 is
the local space $V_{\text{off}}^{\omega_i}$.


\begin{remark}[{On $S$ norm}] \label{rem:Snorm}
In this remark, we show eigenvalues decay
for several choices of $a^{\text{off}}$ and $s^{\text{off}}$.
We consider a target domain $\omega_{t}=[0.4, 0.6]\times [0.4, 0.6]$
 within $D=[0,1]\times [0,1]$. Our space $V_{\text{snapshots}}$ is
generated by solving problems with force terms located outside
$\omega_t$. 
We consider several choices for $a^{\text{off}}$ and $s^{\text{off}}$.
First, we choose 
\[
a^{1,\text{off}}_{\omega_t}=\sum_{i=1}^4 \int_{\omega_t}\kappa(x) |\nabla (\chi_i^0 u)|^2,\ \ 
s^{1,\text{off}}_{\omega_t}= \int_{\omega_t} \kappa(x) |\nabla  u|^2,
\]
where $\chi_i^0$ are bilinear basis functions.
We also consider
\[
a^{2,\text{off}}_{\omega_t}= \int_{\omega_t} \kappa(x) |u|^2,\ \
s^{2,\text{off}}_{\omega_t}= \int_{\omega_t} \kappa(x) |\nabla  u|^2.
\]
We also consider
\[
a^{3,\text{off}}_{\omega_t}= \int_{\omega_t} \kappa(x) |\nabla u|^2,\ \ 
s^{3,\text{off}}_{\omega_{ext}}= \int_{\omega_{ext}} \kappa(x) |\nabla  u|^2,
\]
where $\omega_{ext}=[0.3, 0.7]\times [0.3, 0.7]$. 
The last choice is motivated by  Babuska and Lipton, 2011, \cite{bl11}),
where a larger domain ${\omega_i^+}$, 
${\omega_i}\subset {\omega_i^+}$ (see Figure \ref{coarse_fine_ovs}),
is used
for eigenvalue computations.

We take $\kappa(x)$ to be $10^2$ in $[0.45, 0.55]\times [0.45, 0.55]$
and $1$ elsewhere. In Figure \ref{fig:exampleeigs}, 
we plot the eigenvalues
corresponding to the above choices of $a^{\text{off}}$ and $s^{\text{off}}$.
The snapshot space is generated using unit force terms in the locations
shown in Figure \ref{fig:exampleeigs}, top left.
As we see  in all the cases the eigenvalues decay fast.
The decay of eigenvalues
depends on the choices of $a^{\text{off}}$ and $s^{\text{off}}$, and also
on the choice of $\chi_i^0$.

\begin{figure}[tbp]
\centering
\includegraphics[width=2.4in, height=1.5in]{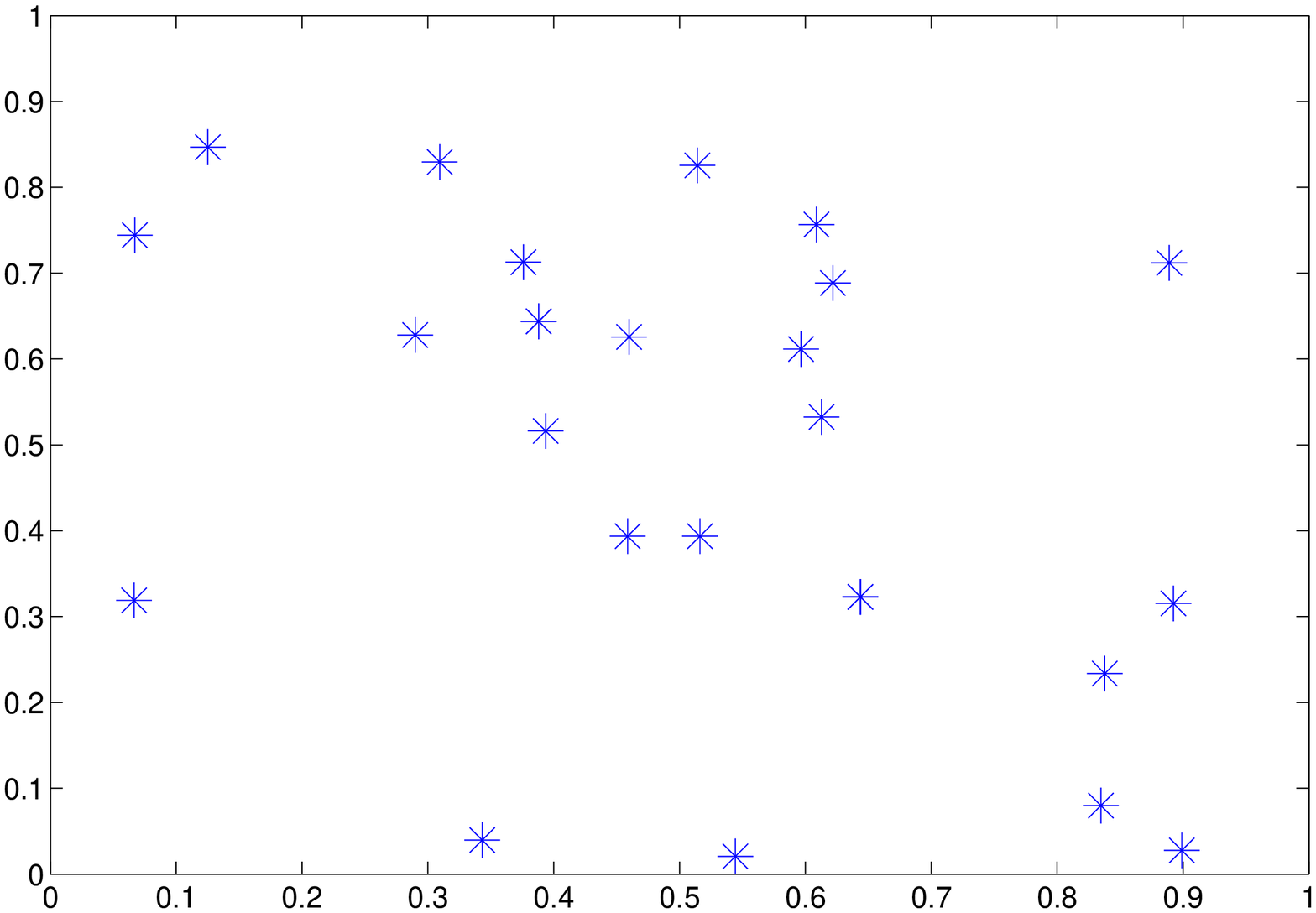}
\includegraphics[width=2.4in, height=1.5in]{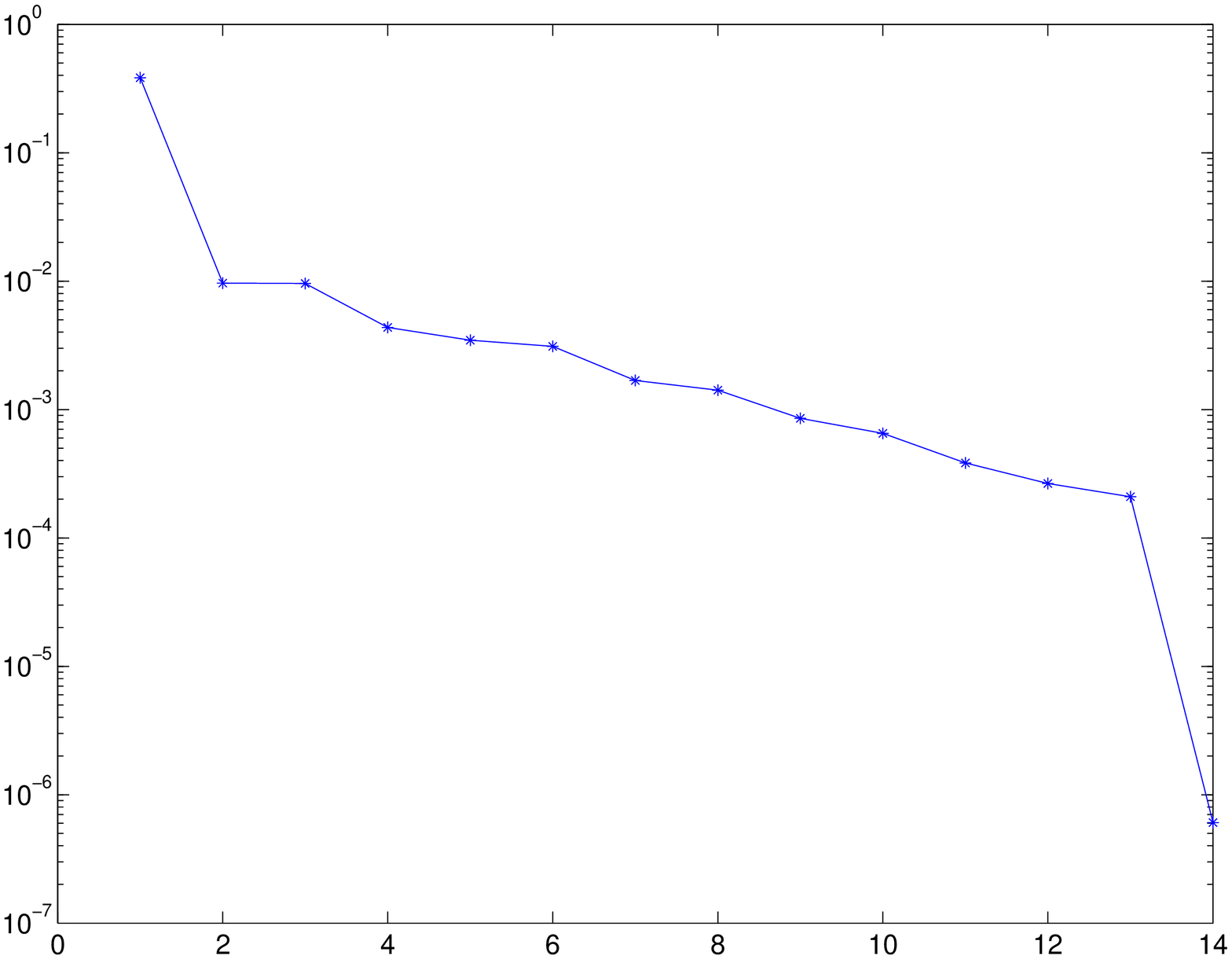}
\includegraphics[width=2.4in, height=1.5in]{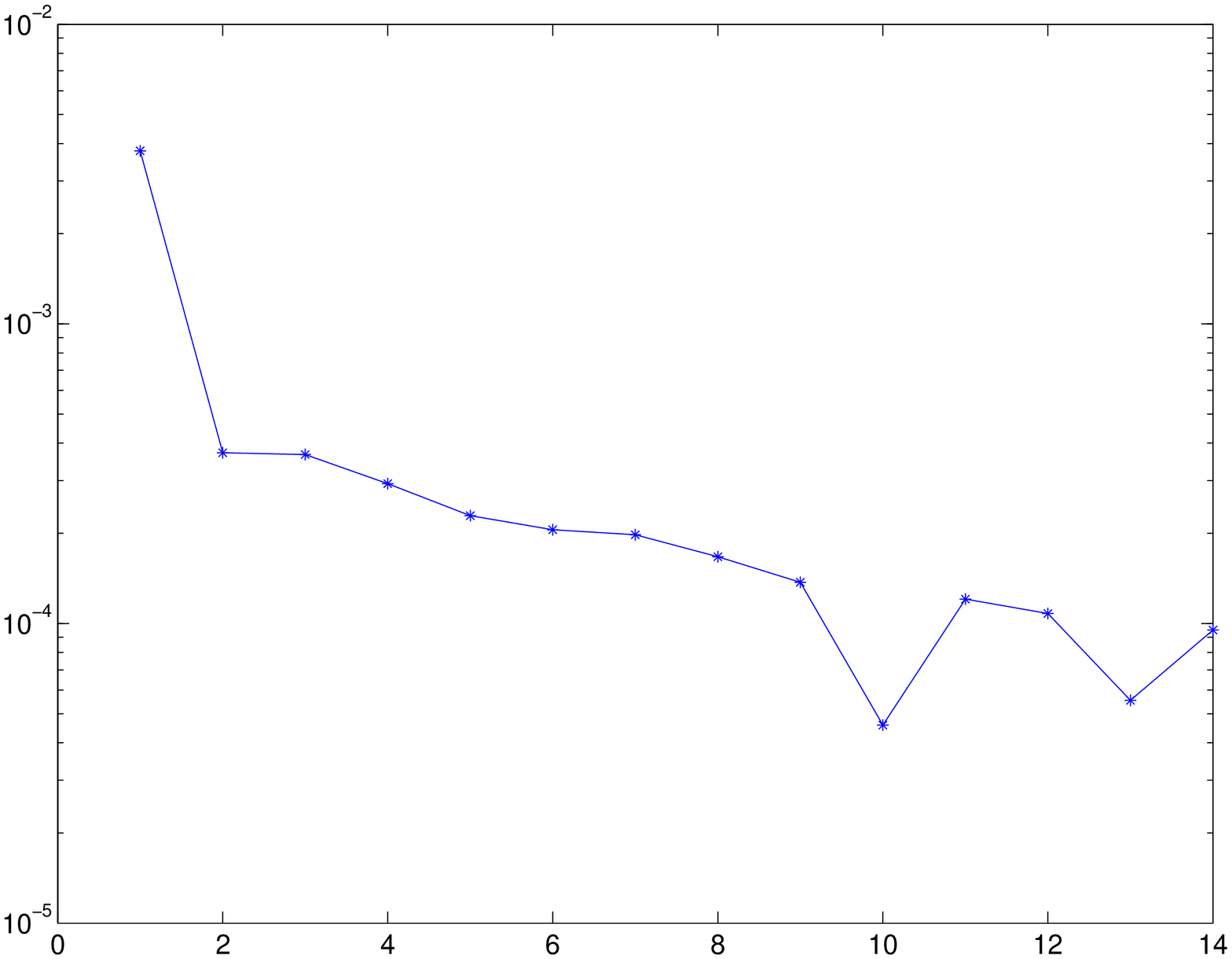}
\includegraphics[width=2.4in, height=1.5in]{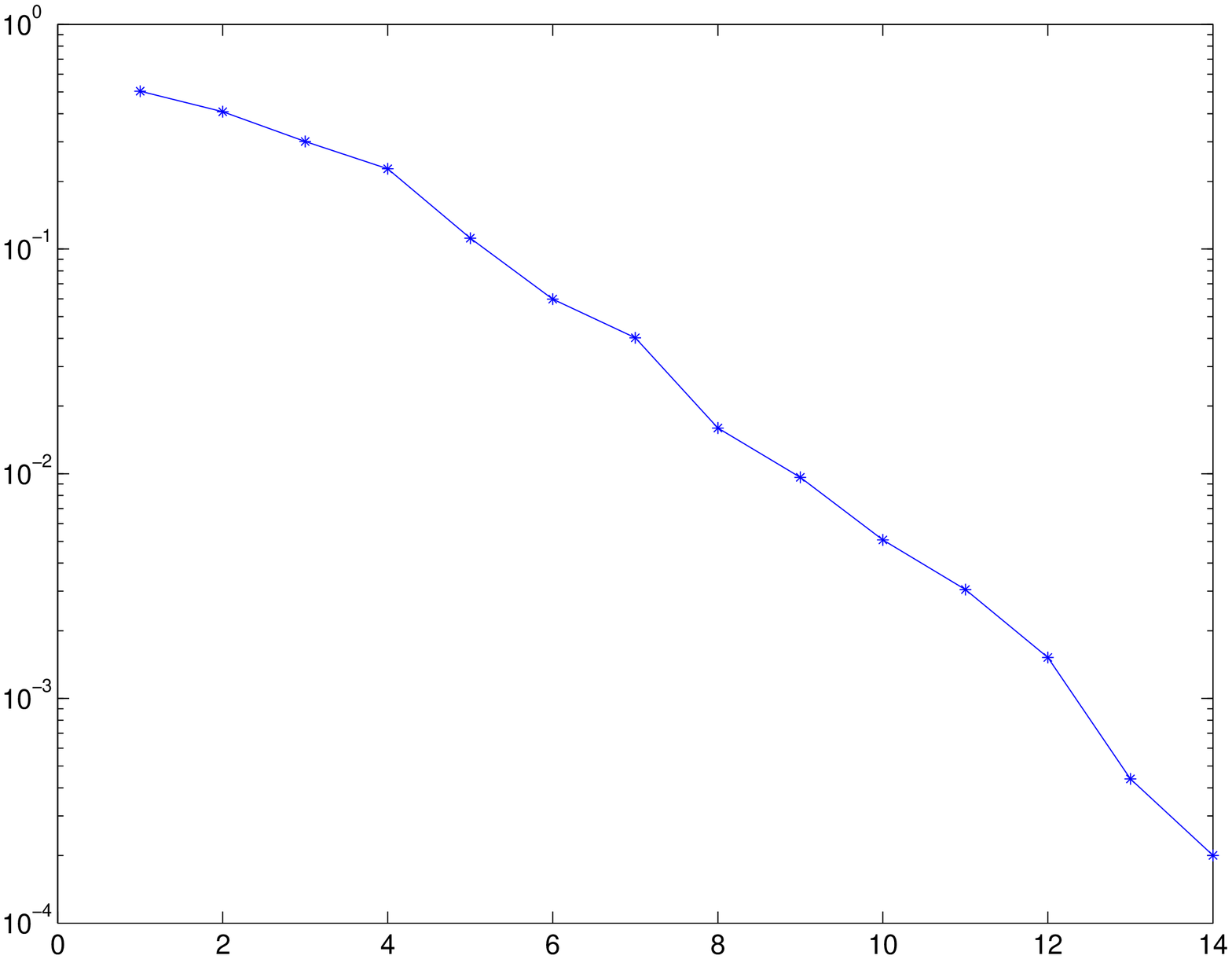}
\caption{Top left: Force locations. Top right: Eigenvalues when using
$a^{1,\text{off}}$
and $s^{1,\text{off}}$. Bottom left: Eigenvalues  when using $a^{2,\text{off}}$
and $s^{2,\text{off}}$. Bottom right: Eigenvalues when using $a^{3,\text{off}}$
and $s^{3,\text{off}}$.}
\label{fig:exampleeigs}
\end{figure}

\end{remark}


\subsection{Step 2. Computing online multiscale basis functions
and their coupling}

At the online stage, for each parameter value, multiscale basis functions
are computed based on each local coarse region. In particular,
for each $\omega_i$ and for each input parameter,
we will formulate a quotient for finding
 a subspace of $V_{\text{on}}^{\omega_i}(\mu)$
where the space will be constructed for each $\mu$ (independent 
of source terms). 
We seek
a subspace $V_{\text{on}}^{\omega_i}(\mu)$ of  $V_{\text{off}}^{\omega_i}$ such that for each
$\psi\in V_{\text{off}}^{\omega_i}$, there exists $\psi_0\in V_{\text{on}}^{\omega_i}(\mu)$ such that
\begin{equation}
\label{eq:eigineq1}
a_{\omega_i}^{\text{on}}(\psi-\psi_0,\psi-\psi_0;\mu)\preceq \delta s_{\omega_i}^{\text{on}}(\psi-\psi_0,\psi-\psi_0;\mu)
\end{equation}
 for some prescribed error tolerance $\delta$ (different
from the one in the offline stage), and the choices
of $a_{\omega_i}^{\text{on}}$ and $ s_{\omega_i}^{\text{on}}$. 
 The corresponding eigenvalue problem is formed
in the space of offline basis functions.
We note that, 
in general, $a_{\omega_i}^{\text{on}}$ and $s_{\omega_i}^{\text{on}}$ contain
partition of unity functions, penalty terms, and other discretization
factors that appear in finite element formulations. In particular,
we assume that
$\kappa_{K}(\chi_i \psi,\chi_i\psi;\mu)\preceq 
a_{\omega_i}^{\text{on}}( \psi,\psi;\mu)$
(where $a_K$ 
corresponds to (\ref{eq:bilinearD}) in $K$).


We note that a choice of partition of unity is important to achieve
smaller dimensional coarse spaces and 
one can look for a choice of optimal partition of unity function.
Because $V_{\text{off}}^{\omega_i}$ is a finite dimensional space, the following
local eigenvalue problem
\begin{equation}
\label{eq:eigon}
A^{\text{on}}\psi_l^{\text{on}}=\lambda_l^{\text{on}} S^{\text{on}} \psi_l^{\text{on}}
\end{equation}
can be used to find the basis functions, where
$A^{\text{on}}$ and  $S^{\text{on}}$ are 
local matrices corresponding to 
$a(\cdot,\cdot;\mu)$ and $s(\cdot,\cdot;\mu)$, 
$A^{\text{on}}=(A^{\text{on}}_{IJ})$ and $S^{\text{on}}=(S^{\text{on}}_{IJ})$ are corresponding
local matrices that are given by 
\[
A^{\text{on}}_{IJ}=a_{\omega_i}^{\text{on}}(\psi_I^{\omega_i},\psi_J^{\omega_i}),\ \
S^{\text{on}}_{IJ}=s_{\omega_i}^{\text{on}}(\psi_I^{\omega_i},\psi_J^{\omega_i}),\ \
\psi_I^{\omega_i}, \psi_J^{\omega_i}\in V_{\text{off}}^{\omega_i}.
\]
 Here, we assume 
that $S^{\text{on}}_{ij}$  is a positive definite matrix.
Note that the matrices $A^{\text{on}}$ and $S^{\text{on}}$ are
computed in the space of offline basis functions.
$V_{\text{on}}^{\omega_i}(\mu)$ is constructed by selecting 
$L_i$ eigenvectors corresponding to the largest
eigenvalues. 
In particular,
the corresponding eigenvectors 
are computed 
by multiplication,  $\sum_{j}\Psi_{ij}^{\text{on}} \psi^{\omega_l}_j$
in each $\omega_l$,  where
$\Psi_{ij}^{\text{on}}$ are coordinates of the vector $\Psi_i^{\text{on}}$.

More precisely, if we reorder the offline
functions using a single index based on the decay of eigenvalues to create the matrices
\[
R_{\text{off}} = \left[ \psi_{1}^{\text{off}}, \ldots, \psi_{M_{\text{off}}}^{\text{off}} \right]
\]
that are defined on the fine grid,
then
\begin{equation*}
 \displaystyle S^{\text{on}}= [s^{\text{on}}_{mn}] = \int_{\omega_i} \kappa(x; \mu) \nabla \psi_m^{\text{off}} \cdot \nabla \psi_n^{\text{off}} = R_{\text{off}}^T {S} R_{\text{off}},
 \end{equation*}
\begin{equation*}
 \displaystyle A^{\text{on}} = [a^{\text{on}}_{mn}] = \int_{\omega_i} {\kappa}(x; \mu) \psi_m^{\text{off}} \psi_n^{\text{off}} = R_{\text{off}}^T {A} R_{\text{off}},
 \end{equation*}
where $S$ and $A$ denote  fine scale matrices
corresponding to the stiffness and mass matrices at $\mu$.
Under certain requirements for the space $V_{\text{off}}^{\omega_i}$,
we can
guarantee that (\ref{eq:eigineq}) holds for all $\psi\in V_{\text{snapshots}}^{\omega_i}$ (and 
not only for all $\psi\in V_{\text{off}}^{\omega_i}$ which is a smaller subspace).



Once multiscale basis functions are constructed, we project the global solution
onto the space of basis functions. One can choose different global 
coupling methods and we present some of them.

{\bf Galerkin coupling.} For a Galerkin formulation, we need conforming
basis functions. We modify $V_{\text{on}}^{\omega_i}$ by multiplying the functions
from this space  with
partition of unity functions. The modified space has the same dimension
and is given by
$Span_{j}(\chi_i \psi_j^{\omega_i,on})$, where
$\psi_j^{\omega_i,on}\in V_{\text{on}}^{\omega_i}(\mu)$
and $\chi_i$ is supported in $\omega_i$.
Then, the Galerkin
approximation can be written as
\[
u^G_{ms}(x;\mu)=\sum_{i,j} c_{j}^i \chi_i(x) \psi_{j}^{\omega_i,\text{on}}(x;\mu).
 \]
If we introduce  
\begin{equation}
\label{eq:G}
V_{\text{on}}^{G}=Span_{i,j}(\chi_i \psi_j^{\omega_i,\text{on}}),
\end{equation}
then Galerkin formulation is given by
\begin{equation}
\label{eq:globalG}
\kappa(u^G_{ms},v;\mu)=(f,v),\ \forall \ v\in V_{\text{on}}^G.
\end{equation}

{\bf Petrov-Galerkin coupling.}
We denote $V_{\text{on}}^{PG}=Span_{i,j}\{ \psi_j^{\omega_i} \}$ and
write the PG approximation of the solution as
\[
u^{PG}_{ms}(x;\mu)=\sum_{i,j} c_{j}^i  \psi_{j}^{\omega_i}(x;\mu).
\]
Then the Petrov-Galerkin formulation is given by
\begin{equation}
\label{eq:globalPG}
\kappa(u^{PG}_{ms},v;\mu)=(f,v),\ \forall \ v\in V_{\text{on}}^G,
\end{equation}
where $V_{\text{on}}^G$ is defined with (\ref{eq:G}).

{\bf Discontinuous Galerkin coupling}

One can also use the discontinuous Galerkin (DG) approach 
(see also \cite{MR1885715,MR2002258,MR2431403})
to couple multiscale basis functions. This may avoid the use of 
the partition of unity functions;
however, a global formulation needs to be chosen carefully.
We have been investigating the use of DG coupling and the detailed
results will be presented elsewhere. Here, we would like to 
briefly mention a general global coupling that can be used. 
The global formulation is given by 
\begin{equation}\label{eq:disc}
{\kappa}^{DG}(u, v) = f(v) \quad \mbox{ for all } 
\quad v=
\{v_{K} \in V_{\text{on}}^K\} ,
\end{equation}
where  
\begin{equation}\label{eq:def:a-h} 
{\kappa}^{DG}(u, v) =\sum_{K} {\kappa}^{DG}_K(u,v)~~~
\mbox{and}~~~~f(v) = \sum_{K} \int_{K} f v_K dx
\end{equation}
for all   
$u=\{u_K\}, v=\{v_K\}$. 
Each local bilinear form ${\kappa}^{DG}_K$ is given as a sum of three 
bilinear forms:
\begin{equation}\label{eq:def:a^hat-i}
 {\kappa}^{DG}_K(u,v) :=  \kappa_K(u, v) + r_K(u, v) + p_K(u, v), 
\end{equation}
where $\kappa_K$ is the bilinear form,  
\begin{equation}\label{eq:def:a-i}
\kappa_K(u, v) := \int_{K}\!\!\kappa_r \nabla u_K \cdot \nabla v_K dx,
\end{equation}
where $\kappa_r$ is the restriction of $\kappa(x)$ in $K$;
the $r_K$ is the symmetric bilinear form, 
\begin{equation}\label{eq:def:s-i}
r_K(u, v) := \sum_ {E\subset \partial K}
 \frac{1}{l_{E}}\int_{E}  \widetilde{\kappa}_{E} \left( \frac{\partial u_K}{\partial n_K} (v_K - v_{K'}) + \frac{\partial v_K}{\partial n_K} (u_{K'} - u_{K}) \right) ds,
\nonumber 
\end{equation}
where $ \widetilde{\kappa}_{E}$ is a weighted average of $\kappa(x)$
near the edge $E$, $l_E$ is the length of the edge $E$,
and $K'$ and $K$ are two coarse-grid elements
sharing the common edge $E$;
and $p_K$ is the penalty bilinear form, 
\begin{equation} \label{eq:def:p-i}
p_K(u, v) := \sum_{E \subset \partial K} 
\frac{1}{l_E} {\delta_E}
\int_{E} 
\widetilde{\kappa}_{E} (u_{K'} - u_K)(v_{K'} - v_K)ds.
\end{equation}
Here  $\delta_E$ is a positive penalty parameter that needs to be selected
and its choice affects the performance of GMsFEM. 
One can choose eigenvalue problems based on DG bilinear forms.
We refer to \cite{DGcoupling} for some results along this direction.

\begin{remark} [{On the convergence of GMsFEM}]
\label{rem:convergence}
Here, we briefly discuss the convergence of the GMsFEM based on Galerking
coupling.
We denote the bilinear form (\ref{eq:bilinearD})
on the whole domain $D$ by $\kappa_D(u,v;\mu)$, and by
$u_{ms}=\sum_{i,j} c_{i,j} \chi_i \psi_j^i$ the GMsFEM interpolant, and 
by $u_0^{\omega_i} =\sum_{j} c_{i,j} \psi_j^i$ an approximation over
the patch $\omega_i$. Then, the error
analysis is the following.
\begin{equation}
\begin{split}
\kappa_D(u-u_{ms},u-u_{ms};\mu)\preceq \kappa_D(u-u_{\text{off}},u-u_{\text{off}}) 
+ \kappa_D(\sum_i \chi_i (u_{\text{off}}-u_0^{\omega_i}),\sum_i \chi_i(u_{\text{off}}-u_0^{\omega_i});\mu)\preceq \\
\kappa_D(u-u_{\text{off}},u-u_{\text{off}};\mu) + \sum_K \kappa_K(\sum_i \chi_i (u_{\text{off}}-u_0^{\omega_i}),\sum_i \chi_i(u_{\text{off}}-u_0^{\omega_i});\mu)\preceq \\
\kappa_D(u-u_{\text{off}},u-u_{\text{off}};\mu) + \sum_i a_{\omega_i}^{\text{on}}( (u_{\text{off}}-u_0^{\omega_i}), (u_{\text{off}}-u_0^{\omega_i});\mu)\preceq \\
\kappa_D(u-u_{\text{off}},u-u_{\text{off}};\mu) + \sum_i \delta s_{\omega_i}^{\text{on}}(u_{\text{off}}-u_0^{\omega_i},(u_{\text{off}}-u_0^{\omega_i});\mu)\preceq\\
\kappa_D(u-u_{\text{off}},u-u_{\text{off}};\mu) + \sum_i \delta s_{\omega_i}^{\text{on}}(u,u;\mu) \preceq \kappa_D(u-u_{\text{off}},u-u_{\text{off}};\mu) + \delta s_D^{\text{on}}(u_{\text{off}},u_{\text{off}};\mu).
\end{split}
\end{equation}
Here, we used the fact that 
$\kappa_{K}(\chi_i \psi,\chi_i\psi;\mu)\preceq 
 a_{\omega_i}^{\text{on}}(\chi_i \psi,\chi_i \psi;\mu)$,
which is an assumption on the choice of $ a_{\omega_i}^{\text{on}}$,
the inequality (\ref{eq:eigineq}), and replaced the sum over $K$
by the sum over larger regions $\omega_i$. We note that the assumption
 $\kappa_{K}(\chi_i \psi,\chi_i\psi;\mu)\preceq 
 a_{\omega_i}^{\text{on}}(\chi_i \psi,\chi_i \psi;\mu)$
can be easily satisified with an appropriate choice of
$a_{\omega_i}^{\text{on}}$ in continuous Galerkin framework (see \cite{egw10}).
One can choose $s_{\omega_i}^{\text{on}}(\chi_i \psi,\chi_i \psi;\mu)$
 based on finite element formulation such that
$s_{\omega_i}^{\text{on}}(\chi_i \psi,\chi_i \psi;\mu)=\kappa_{\omega_i}(\chi_i \psi,\chi_i \psi;\mu)$.
For discontinuous Galerkin formulation (with no partition of unity
$\chi_i$), we refer to \cite{DGcoupling} for the appropriate choice of 
$a_{\omega_i}^{\text{on}}$.
The hidden constants depend on the number of neighboring 
element of each coarse block and the  constants in (\ref{eq:eigineq}) and
(\ref{eq:eigineq1}). Moreover, there is an irreducable error
$\kappa_D(u-u_{\text{off}},u-u_{\text{off}};\mu)$.
We note that one can obtain sharper estimates using a bootstrap argument
in parameter-independent case (see \cite{egw10, eglp13})
under some assumptions
 and show that the error decreases
as the coarse-mesh size decreases.

\end{remark}

\subsection{Iterative solvers - online correction of fine-grid solution}

In the previous approach, the coarse spaces are designed
 to achieve a desired accuracy.
One can also iterate (on residual and/or the basis 
functions) for a given source term and converge to the true solution
without increasing the dimension of the coarse space. 
Both approaches have their application
fields and allows computing the fine-grid solution with an increasing
accuracy.
Next, we briefly describe the solution procedure based on 
the concept of two-level iterative methods.

We denote by  $\{D_i'\}_{i=1}^M$ the overlapping decomposition obtained from the original  non-overlapping decomposition $\{ D_i\}_{i=1}^M$ by enlarging each subdomain $D_i$ to 
\begin{equation}\label{eq:def:overlapping-subdomains}
D'_i=D_i\cup \{x\in D, \mbox{dist}(x,D_i)<\delta_i\}, \quad i=1,\dots,M,
\end{equation}
where $\mbox{dist}$ is some distance function and let $V^h_0(D_i')$ be the set of finite element functions with support in $D_i'$ and zero trace on the boundary $\partial D_i'$. We also denote by $R_i^T:V^h_0(D_i')\to V^h$ the 
extension by zero operator.

With the help of $(u-u_0)$, we can
 correct the coarse-grid solution. A number of approaches can be used. 
For instance, we write the solution in the form
\[
u=u_0 + \sum_i \widehat{\chi}_i v_i,
\]
where $v_i$ are defined in $\omega_i$
(though it can be taken to be supported in a different domain)
and $\widehat{\chi}_i$ is 
a partition of unity.
Suppose that $v_i$ has zero trace on $\partial \omega_i$. 
We can solve the local problems
\[
L(v_i)=f-L_\mu(u_0),
\]
with zero boundary condition.
Other correction schemes can be implemented to correct 
the coarse-grid solution. For example, we 
can use the traces of $u_0$ to correct the solution 
in each subdomain in a consecutive fashion.

When the bilinear form 
is symmetric and positive definite, we can 
consider additive two-level domain 
decomposition  methods  to find the solution of the
fine-grid finite element   problem 
\begin{equation}
\label{eq:fine-system}
\kappa(u,v;\mu)=f(v),\quad \mbox{for all}\  v\in V^h,
\end{equation}
where $V^h$ is the fine-grid finite element space 
of piecewise linear polynomials. The matrix  
of this linear systems is written as
\[
S(\mu) u(\mu) =b.
\]
Here, $S$ is the stiffness matrix associated to the bilinear 
form $\kappa$ and $b$  such that $v^Tb=f(v)$ for all $v\in V^h$.
We can solve the fine-scale linear system 
iteratively with the preconditioned 
conjugate gradient (PCG) method (or any other 
Krylov type method for a non-positive problem).
Any other suitable iterative scheme can be used as well. 
We introduce the two level additive preconditioner of the form 
\begin{equation}\label{eq:def:Binv}
B^{-1}(\mu)=R_{0,on}^T(\mu) \widehat{S}_{0}^{-1}(\mu) R_{0,\text{on}}(\mu)+\sum_{i=1}^M R_i^T S_i^{-1}(\mu) R_i, 
\end{equation}
where the local matrices are defined by 
\begin{equation}\label{eq:def:A^Di}
v^TS_i(\mu)w=a(v,w;\mu) \quad \mbox{ for all } v,w\in V^h_0(D'_i).
\end{equation}
 The coarse projection matrix $R_{0,\text{on}}^T$ is defined by
$R_{0,\text{on}}^T=R_{0,\text{on}}^T(\mu)=[\Phi_1,\dots,\Phi_{N_v}]$ and the online 
coarse matrix 
$\widehat{S}_0(\mu)=R_{0, \text{on}}(\mu)S(\mu)R_{0,\text{on}}^T(\mu)$. 
The columns $\Phi_i$'s are fine-grid coordinate 
vectors corresponding 
to the  basis functions, e.g., in the Galerkin formulations 
they correspond to the basis functions $\{\chi_i(x) \psi_{j}^{\omega_i,\text{on}}(x;\mu)\}$.  
See \cite{tw,Tarekbook} 
and references therein for more details on various domain 
decomposition methods.

The application of the preconditioner involves
 solving local problems in each iteration. In 
domain decomposition methods, our main goal is to reduce
the number of iterations in the iterative procedure. 
It is well known that a coarse solve needs to be added to 
the one level preconditioner in order to construct 
robust methods.  The appropriate construction of the coarse 
space $V_0$  plays a key role in 
obtaining robust iterative domain decomposition methods.
Our methods provide an inexpensive coarse solves 
and efficient iterative solvers for general parameter-dependent 
problems. This will be discussed in the next sections.

\section{Case studies and relation to existing methods.  Discussions and applications}

In this section, we illustrate basic concepts via some specific examples.
We use existing methods in the literature for multiscale
problems and show how these methods can be put under the general framework
of the GMsFEM and show some numerical results.

\subsection{Case with no parameter}
\label{sec:noinput}

In this section, we write a method proposed in \cite{ge09_2} as
a special case of GMsFEM.
First, we consider a case with no parameter, i.e.,
\[
L(u)=f,\ \text{or corresponding} \ \kappa(u,v)=f(v).
\]
In this case, offline multiscale basis functions are used
for online simulations due to the absence of the parameter.
Next, we discuss the construction of
multiscale basis functions
(see \cite{ge09_2, eglw11,Review}).

The construction of the offline space starts with the snapshot space.
We choose
the snapshot space as all fine-grid functions within a coarse region.

For the construction of the offline space, we can choose 
\begin{equation}
\label{eq:defa}
a_{\omega_i}^{\text{off}}(\psi,\psi)=\sum_{k,\omega_k\bigcap\omega_i \not= 0} \kappa_{\omega_k}(\chi_k \psi,\chi_k \psi),
\end{equation}
where $\kappa_{\omega_i}(\cdot,\cdot)$ is the restriction of $\kappa(\cdot,\cdot)$
in $\omega_i$ and
$\chi_k$ is the partition of unity corresponding to $\omega_k$. 
For example, for the elliptic equation (see (\ref{eq:ellipcase})),
\[
\kappa_{\omega_i}(\psi,\psi)=\int_{\omega_i} \kappa \nabla \psi \cdot \nabla \psi.
\]
As for $s_{\omega_i}^{\text{off}}$, we select
\[
s_{\omega_i}^{\text{off}}(\psi,\psi)= \kappa_{\omega_i}( \psi, \psi).
\]
The corresponding eigenvalue problem can be explicitly written as before.
In our numerical implementation, we choose
$a_{\omega_i}^{\text{off}}(\psi,\psi)=\sum_{k,\omega_k\bigcap\omega_i \not= 0} \int_{\omega_i} \kappa |\nabla \chi_k|^2 |\psi|^2$ (see \cite{egw10} that shows that this is not smaller
than (\ref{eq:defa}) in the space of local $\kappa$-harmonic functions).

In the first example, we do not construct an online space and the 
offline space is used in GMsFEM, $V_{\text{on}}^{\omega_i}=V_{\text{off}}^{\omega_i}$.
The basis functions are constructed by
multiplying 
the eigenvectors corresponding to the dominant eigenvalues
by partition of unity functions, see (\ref{eq:G}). The Galerkin
coupling of these basis functions are performed based on 
(\ref{eq:globalG}). Note that the stiffness is pre-computed
in the offline stage and there is no need for any stiffness matrix
computation in the online stage.

We point out that the choice of initial partition of unity basis functions
$\chi_i$, are important in reducing the number of very large
eigenvalues. 
We note that the dimension of the coarse space depends on the choice
of $\chi_i$ and, thus, it is important to have a good
choice of  $\chi_i$.
The essential ingredient in designing  them is
to guarantee that there are fewer large
eigenvalues, and thus the coarse space dimension is small.
With an initial choice of multiscale basis functions $\chi_i$
that contain many localizable small-scale  features of the solution,
one can reduce the dimension of the resulting coarse space.

\begin{figure}[htb]
\centering
{\includegraphics[width=6cm, height=5cm]{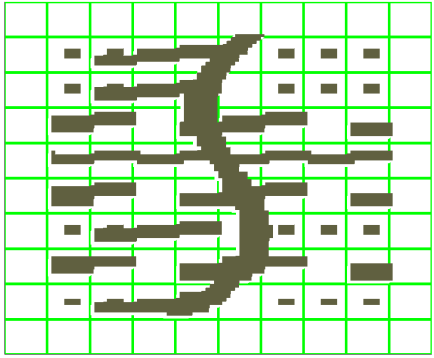}}
\caption{Coefficient  $\kappa(x)$. 
The dark region has high conductivity $\eta$
and the white background has conductivity 1. Green lines show
the coarse grid.}
\label{fig:newperm2} 
\end{figure}

Next, we briefly discuss a few numerical examples
(see \cite{ge09_2} for more discussions). 
We present a numerical result
for the  coarse-scale approximation
and for the two-level additive preconditioner
(\ref{eq:def:Binv}) with the 
local spectral multiscale coarse spaces as discussed above.
The equation 
$-\mbox{div}(\kappa\nabla u)=1$ is solved with boundary conditions 
$u=x+y$ on $\partial D$.
For the coarse-scale approximation, we  vary the dimension
of the coarse spaces by adding additional basis
functions corresponding to the largest eigenvalues.
We investigate the convergence rate, while
for preconditioning results,
we will investigate the behavior of the condition number
as we increase the contrast for various choices of coarse spaces. 
The domain 
$D=[0,1]\times [0,1]$ is divided into  $10\times 10$ equal square 
subdomains. 
Inside each subdomain we use a fine-scale triangulation, where 
triangular elements constructed from $10\times 10$ squares are
used. We  consider the scalar coefficient 
$\kappa(x)$ depicted in Figure \ref{fig:newperm2}
that corresponds to a background one and high 
conductivity channels
 and inclusions. 

We test the accuracy of GMsFEMs when coarse spaces include eigenvectors
corresponding to the large eigenvalues.
We implement GMsFEM by choosing the initial partition of unity functions
to consist
of multiscale functions with linear boundary conditions (MS)
(see (\ref{eq:def:msfem}).
We use the following notation.
GMsFEM$+0$ refers to the GMsFEM where the coarse space
includes all eigenvectors that correspond to eigenvalues which
are asymptotically unbounded as the contrast increases, i.e.,
these eigenvalues increase as we increase the contrast.
One of these eigenvectors corresponds to a constant function
in the coarse block.
GMsFEM$+n$ refers to the GMsFEM where in addition to
eigenvectors that correspond to asymptotically unbounded eigenvalues,
we also add $n$
eigenvectors corresponding to the next $n$ eigenvalues. 

In previous studies \cite{ge09_1, ge09_3, egw10}, 
we discussed how the number
of these asymptotically unbounded eigenvalues depends
on the number of inclusions and channels. In particular,
we showed that if there are $n$ inclusions (isolated regions
with high conductivity) and $m$ channels (isolated high-conductivity
regions connecting boundaries of a coarse grid), then the number
of asymptotically unbounded eigenvalues is $n+m$ when
standard bilinear partition of unity function, $\chi_i^0$,  used.
However, if 
the partition of unity $\chi_k$ is chosen as multiscale finite element
basis functions (\cite{hw97}) defined by
\begin{eqnarray}
\mbox{div}(\kappa\nabla\chi_i^{ms})=0\ \ \mbox{in }K\in \omega_i,\quad
\chi_i^{ms}=\chi_i^0\ \ \mbox{in }\partial K,\ \ \forall\  K\in \omega_i,
\label{eq:def:msfem}
\end{eqnarray}
then the number of asymptotically unbounded eigenvalues is $m$.
We can also use energy minimizing basis functions that are defined
(see \cite{XuZikatanov})
as
\begin{equation}\label{eq:energyminization}
\min\sum_{i}
\int_{\omega_i}\kappa |\nabla \chi_{i}^{emf}|^2
\end{equation}
subject to $\sum_{i} \chi_i^{emf}=1$
 with $\mbox{Supp}(\chi_i)\subset \omega_i$, 
$i=1,\dots,N_v$, to achieve even smaller dimensional coarse spaces.

In all numerical results, the errors are measured in
 the energy norm ($|\cdot|_A^2$), 
 $H^1$ norm ( $|\cdot|_{H^1}^2$), and
$L^2$-weighted  norm ($|\cdot|_{L^2}^2$) respectively.
We present the convergence as we increase 
the number of additional eigenvectors. 
In Table  \ref{tab:newperm2mstilde-10to6},
we present the numerical results when the initial 
partition of unity 
consists of multiscale basis functions with linear boundary conditions
for the contrast  $\eta=10^6$.   We note that the convergence 
is robust with respect
to the contrast and the error reduces. The error is
proportional to the largest eigenvalue ($\Lambda^*$) whose
eigenvector is not included in the coarse space as one 
can observe from the table (correlation coefficient between
$\Lambda^*$ and the energy error is $0.99$).
We observe that the errors are smaller compared to those obtained using 
MsFEM with
piecewise linear initial conditions.

\begin{table}
\centering\small 
\begin{tabular}{|l|r|r|r|r|r|r|r|}\hline
$H=1/10$ & $|\cdot|_A^2$& $|\cdot|_{H^1}^2$ & 
$|\cdot|_{L^2}^2$ & $\Lambda^*$ \\\hline
{GMsFEM}+0 (153) & 14.3  &14.3  & 6.17e-02   & 0.0704\\
{GMsFEM}+1 (234) & 5.82 & 5.82 &  1.02e-02    &0.0117\\
{GMsFEM}+2 (315) &  5.33& 5.33 & 8.60e-03    & 0.0071\\
{GMsFEM}+3 (396) &  4.64& 4.64 & 6.53e-03    &0.0043\\
{GMsFEM}+4 (477) & 4.01 & 4.01 & 4.80e-03     &0.0032\\\hline\hline
\end{tabular}
\caption{
Convergence results (in \%) for GMsFEM with MS with  increasing dimension
of the coarse space. 
Here, 
$\eta=10^6$. The initial coarse space is spanned by
multiscale basis functions
with piecewise linear boundary conditions ($\chi^{ms}$).
The coefficient is depicted in Figure 
\ref{fig:newperm2}.}
\label{tab:newperm2mstilde-10to6}
\end{table}

Next, we present numerical results when the snapshot space consists
of harmonic functions in $\omega_i$. 
More precisely, for each fine-grid function, $\delta_l^h(x)$,
which is defined by 
$\delta_l^h(x)=\delta_{l,k},\,\forall l,k\in \textsl{J}_{h}(\omega_i$, where $\textsl{J}_{h}(\omega_i^{+})$ denotes the fine-grid boundary node on $\partial\omega_i^{+}$,
we
solve
\[
-div(\kappa (x)\nabla  \psi_{l}^{\text{snap}})=0\ \ \text{in} \ \omega_i
\]
subject to boundary condition, $ \psi_{l}^{\text{snap}}=\delta_l^h(x)$.
We use bilinear partition of unity functions for the mass matrix. The numerical
results are presented in Table \ref{tab:newperm2mstilde-10to6_harmonic}
and we observe slightly worse results compared to
 Table \ref{tab:newperm2mstilde-10to6}; however, the errors are comparable
for dimensions of order $400$.
We have observed similar results if larger regions are used for computing 
snapshot functions (see Remark \ref{rem:Snorm}).
We report the results on oversampling techniques in \cite{eglp13}

\begin{table}
\centering\small 
\begin{tabular}{|l|r|r|r|r|r|r|r|}\hline
Coarse dim & $|\cdot|_A^2$& $|\cdot|_{L^2}^2$  \\\hline
dim=202 & 13.7  & 0.2   \\
dim=364 & 3.24 &   1.44e-3   \\
dim=607 &  0.021&  1.0e-4   \\
dim=850 & 0.164 &  5.6e-5    \\\hline\hline
\end{tabular}
\caption{
Convergence results (in \%) for GMsFEM with MS with  increasing dimension
of the coarse space. 
Here, 
$\eta=10^6$. The initial coarse space is spanned by
bilinear basis functions and harmonic functions are used as a space of snapshots.
The coefficient is depicted in Figure 
\ref{fig:newperm2}.}
\label{tab:newperm2mstilde-10to6_harmonic}
\end{table}

Next, we present the results for a two-level preconditioner.
We implement  a two-level additive preconditioner with  
the following 
coarse spaces: 
multiscale functions with linear boundary conditions (MS);
energy minimizing functions (EMF);
spectral coarse spaces using
piecewise linear
partition of unity functions as an initial space (GMsFEM  with Lin);
spectral coarse spaces 
where 
 multiscale finite
element basis functions with linear boundary conditions
($\chi^{ms}$)
are used as an initial partition of unity
(GMsFEM  with MS);
spectral coarse spaces with $\widetilde{\kappa}$ 
where
energy minimizing
basis functions ($\chi^{emf}$) are used as an initial 
partition of unity
(GMsFEM with EMF).
In Table \ref{tab:newperm2}, we show the number of PCG iterations
and estimated condition numbers. We also show
the dimensions of the coarse spaces. 
Note that
the standard coarse space with one basis per coarse node 
has the dimension $81\times 81$. The smallest dimension can be
achieved by
using energy minimizing basis functions as an 
initial partition of unity.
We observe that the number of iterations does not change as
the contrast increases when spectral coarse spaces are used.
This indicates that the preconditioner is optimal.
On the other hand, when using multiscale basis functions 
(one basis per coarse node),
the condition number of the preconditioned matrix increases
as the contrast increases.

\begin{table}[htb]
\centering
\footnotesize
\begin{tabular}{|l|r|r|r|r|r|r|}\hline
$\eta$ & MS& EMF   &GMsFEM with Lin
& GMsFEM with MS
 &   GMsFEM with EMF  \\\hline
 $10^3$    &83(2.71e+002)   &69(1.43e+002)   &31(8.60e+000)   &31(9.34e+000)   &32(9.78e+000)  \\ 
  $10^5$    &130(2.65e+004)   &74(1.29e+004)   &33(8.85e+000)   &33(9.72e+000)   &34(1.02e+001)  \\ 
  $10^7$    &189(2.65e+006)   &109(1.29e+006)   &34(8.85e+000)   &35(9.60e+000)   &37(1.02e+001)  \\ \hline\hline
 Dim &  81& 81 &165 &113  &113\\
\hline\hline 
\end{tabular}
\caption{Number of iterations until convergence 
and estimated condition number for  the PCG  and different values of 
the contrast $\eta$ with the coefficient depicted in  
Figure \ref{fig:newperm2}. We set 
the tolerance to 1e-10. Here $H=1/10$ 
with $h=1/100$.}
\label{tab:newperm2}
\end{table}

\subsection{Elliptic equation with input parameter}
\label{sec:RB}

Here, we present a method proposed in \cite{egt11} as a special case of 
GMsFEM.
We consider a parameter-dependent elliptic equation (see (\ref{eq:ellipcase}))
\begin{equation}\label{eq:general-parametric-problem}
L_\mu(u)=f,\ \text{or corresponding} \ \kappa(u,v;\mu)=f(v).
\end{equation}
As before, we start the computation with the snapshot space
consisting of local fine-grid functions
and compute offline basis functions. In this example, we
will construct offline multiscale
spaces for some selected values of $\mu$, $\mu_i$ ($i=1,...,N_{rb}$),
where $N_{rb}$ is the number of selected values of $\mu$ used
in constructing multiscale basis functions. 
These values of $\mu$ are selected via an inexpensive RB
 procedure
\cite{egt11}. We
briefly describe the offline space construction.
For each selected $\mu_j$ (via RB procedure \cite{egt11}), we choose the offline space 
\[
a_{\omega_i,j}^{\text{off}}(\psi,\psi)=\int_{\omega_i} \kappa(x;\mu_j)  \psi \psi.
\]
As for $s_{\omega_i}^{\text{off}}$, we select
\[
s_{\omega_i,j}^{\text{off}}(\psi,\psi)=\kappa_{\omega_i}(\psi,\psi;\mu_j).
\]
Then, the selected dominant eigenvectors are orthogonalized
with respect to $H^1$ inner product. 


At the online stage, for each parameter value, multiscale basis functions
are computed based on the offline space. In particular,
for each $\omega_i$ and for each input parameter,
 we formulate a quotient to find
 a subspace of $V_{\text{on}}^{\omega_i}(\mu)$,
where the space will be constructed for each $\mu$. 
For the construction of the online space, we choose 
\[
a_{\omega_i}^{\text{on}}(\psi,\psi;\mu)=\int_{\omega_i} \kappa(x;\mu)  \psi^2.
\]
For the bilinear form for $s$, we choose
\[
s_{\omega_i}^{\text{on}}(\psi,\psi;\mu)=\kappa_{\omega_i}(\psi,\psi;\mu).
\]

In this case, the online space is a subspace of the offline
space and computed by solving an eigenvalue problem for a given value
of the parameter $\mu$. The online space is computed by solving
an eigenvalue problem in $\omega_i$ using $V_{\text{off}}^{\omega_i}$,
see (\ref{eq:eigon}). Using dominant eigenvectors, we form
a coarse space as in (\ref{eq:G}) and solve the global coupled system
following (\ref{eq:globalG}). For the numerical
example, we will consider the coefficient that has an affine
representation (see (\ref{hyp})). The online computational
cost of assembling the stiffness matrix involves summing
$Q$ pre-computed matrices corresponding to coarse-grid systems.
We point out that the choice of initial partition of unity functions,
$\chi_i$, are important in reducing the number of very large
eigenvalues (we refer to \cite{egt11} for further discussions).

We present a numerical example for 
$-\mbox{div}(\kappa(x;\mu)\nabla u)=1$ which is solved with boundary conditions 
$u=x+y$ on $\partial D$.
We take 
$D=[0,1]\times [0,1]$ that is
divided into  $10 \times 10$ equal square subdomains. 
As in Section \ref{sec:noinput},
in each subdomain we use a fine-scale triangulation, where 
triangular elements constructed from $10\times 10$ squares are
used.

We consider 
a permeability field which is the sum of four permeability fields each of
which contains inclusions such that their sum gives
several
channelized permeability field scenarios. 
The permeability field is described by
\begin{equation}\label{eq:kappaxmunumerics1}
\kappa(x;\mu) := \mu_1\kappa_1(x) + \mu_2\kappa_2(x)+\mu_3\kappa_3(x)+\mu_4\kappa_4(x).
\end{equation}
There are several distinct features in this family of 
conductivity fields which include inclusions and the channels
that are obtained by choosing $\mu_1=\mu_2=1/2$ or
$\mu_3=\mu_4=1/2$.
There exists  no single value 
of $\mu$ that has all the features. Furthermore, we will use a 
trial set for the reduced basis algorithm that does not include 
$\mu_i=1/2$ $(i=1,2,3,4)$. 

%
%

\begin{figure}[htbb]
\center
\begin{tabular}{ccc}
\includegraphics[width=7cm]{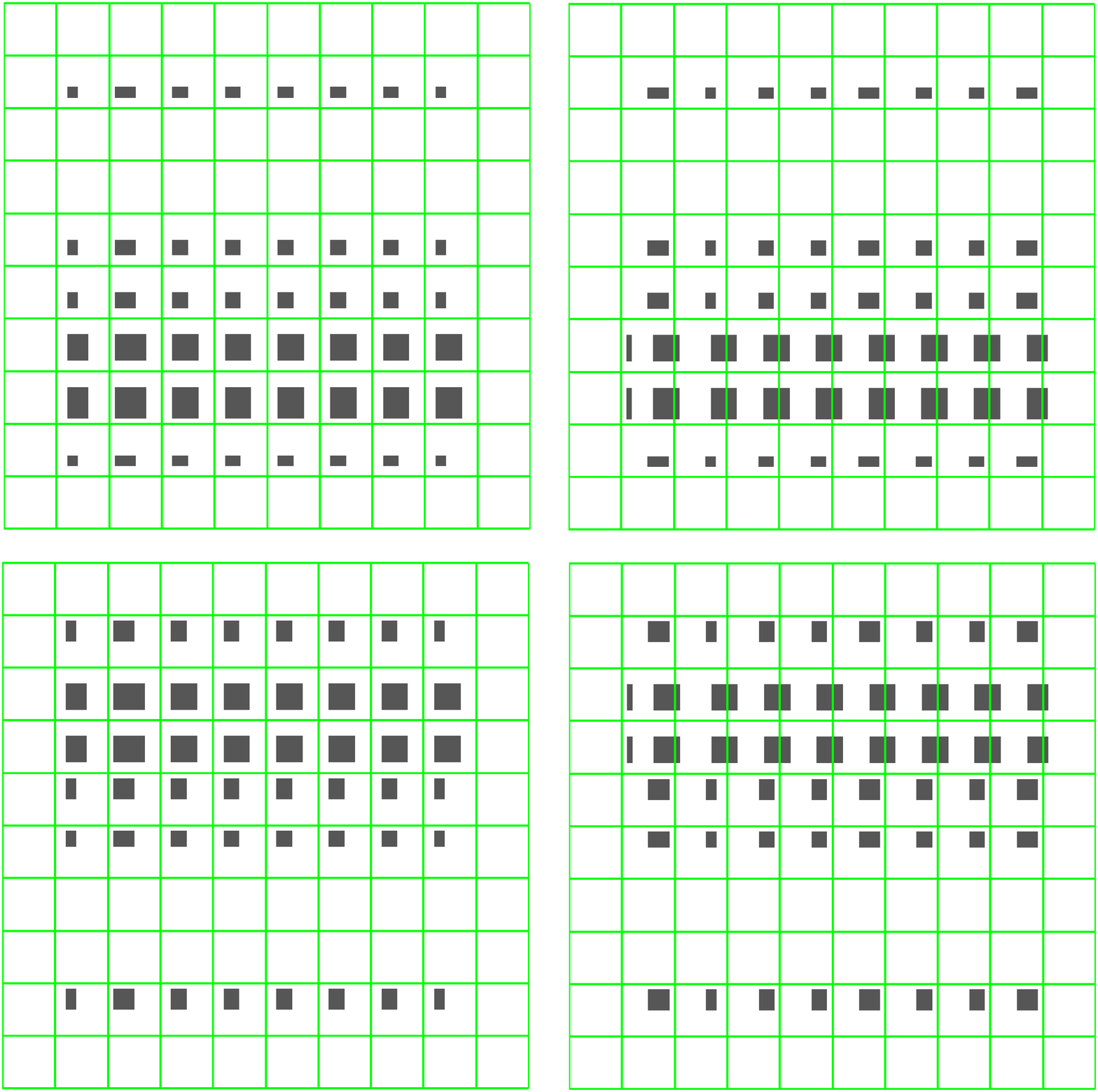}  
\end{tabular}
\caption{From left to right and top to bottom: $\kappa_1$, $\kappa_2$, $\kappa_3$, $\kappa_4$.}
\label{illus}
\end{figure}


As there are several distinct spatial fields in the space of conductivities,
we will choose multiple functions  in our 
reduced basis. 
In the case of insufficient number of samples in the offline space,
we observe that the online permeability field does not contain
appropriate features and this affects the convergence rate.
We observe in 
Table~\ref{tab:newperm1ms-10to6} 
that we indeed need $N_{rb} \geq 4$ to capture all the details 
of the solution. 
In this table, we compare the errors obtained by GMsFEM
with a different number of online basis functions. 
We observe convergence with respect 
to the number of local eigenvectors when $N_{rb}$ increases.
We note that in these computations, we also have an error
associated with the fact that the offline space is not sufficiently
large and thus the error decay is slow as we increase the number
of basis functions.
We have also computed weighted $L_2$ error
which shows a similar trend and the $L_2$ errors are generally
 much smaller.



\begin{table}
\centering\small 
\begin{tabular}{|l|r|r|r|r|r|r|r|r|r}\hline
$H=1/10$         & $N_{rb}=2$ & $N_{rb}=3$ & $N_{rb}=4$  \\\hline
GMsFEM+0             &  $45.5(478)$ &  $40.3(565)$ & $9.2(588)$ \\
GMsFEM+1             &  $39.4(599)$ &  $27.5(686)$ & $5.3(709)$ \\
GMsFEM+2             &  $38.4(720)$ &  $26.8(807)$ & $5.1(830)$\\
GMsFEM+3             &   $36.2(841)$& $26.2(928)$   & $4.9(951)$\\\hline\hline 
\end{tabular}
\caption{Convergence results (energy norm in $\%$ and space dimension) for GMsFEM with the increasing dimension
of the coarse space. 
Here, $h=0.01$, 
$\eta=10^6$, and $\mu_1 =\mu_2=\mu_3=\mu_4=1/2$ (error with MsFEM $96.77\%$).}
\label{tab:newperm1ms-10to6}
\end{table}



\subsection{Anisotropic flows in parameter-dependent media}

\begin{figure}[htbb]
\center
\begin{tabular}{ccc}
\includegraphics[width=11cm]{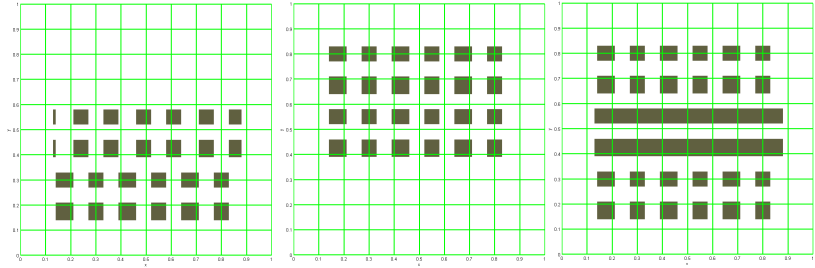}  
\end{tabular}
\caption{From left to right: $\mu=0$, $\mu=1$ and $\mu=1/2$.}
\label{illus-onedim}
\end{figure}

In this example, we apply GMsFEM to anisotropic flows by
considering the elliptic problem with tensor coefficients 
\[
\kappa(x,\mu)=\left(\begin{array}{cc} \kappa_{11}(x,\mu)&0\\ 0&1\end{array}\right)\\
\]
where the $\kappa_{1,1}$ coefficient  is described by
\begin{equation}\label{eq:kappaxmunumerics2}
\kappa_{1,1}(x;\mu) := (1-\mu) \kappa_0(x) + \mu \kappa_1(x).
\end{equation}
We consider an example from where the first component of
the permeability has $3$ distinct different features in $\kappa(x;\mu)$: 
inclusions (left), channels (middle), and
shifted inclusions (right); see Figure~\ref{illus-onedim}. 
The permeability field is described by
\begin{equation}\label{eq:kappaxmunumerics}
\kappa(x;\mu) := (1-\mu) \kappa_0(x) + \mu \kappa_1(x).
\end{equation}
We can represent $3$ distinct different features in $\kappa(x;\mu)$: 
inclusions (left), channels (middle), and
shifted inclusions (right), see Figure~\ref{illus-onedim}. 
There exists  no single value 
of $\mu$ that has all the features. Furthermore, we will use a 
trial set for the reduced basis algorithm that does not include 
$\mu=0.5$.

We will use a 
trial set for the reduced basis algorithm that does not include 
$\mu=0.5$. The trial set is chosen as in the case of isotropic 
case using a greedy algorithm.
The important features in these permeability fields 
characterize the preferred directions of conductivity.

We note that the coarse space has a structure that differs from
the case of isotropic coefficients. In particular, 
the coarse space has a larger dimension and contains
all functions that are constant along $x_1$ direction
in the examples under consideration \cite{eglmr11}.

First, we present numerical results for the
two-level domain decomposition solvers.
In Table \ref{tab:newperm2-0ANI}, 
we present the results for a two-level preconditioner.
We implement  a two-level additive preconditioner with  
the coarse spaces introduced earlier: 
multiscale functions with linear boundary conditions (MS);
spectral coarse spaces where piecewise linear
partition of unity functions are used as an initial 
partition of unity
(GMsFEM with Lin); and
spectral coarse spaces where multiscale finite
element basis functions with linear boundary conditions
($\chi^{ms}$)
are used as an initial partition of unity
(GMsFEM with MS).

From this table, we see that the number of PCG iterations
and estimated condition numbers do not depend on the contrast
when spectral basis functions are used. 
We use $N_{rb}=3$ because we need at least three features to represent
the permeability field as in the earlier example. We 
observed from our previous study that if $N_{rb}<3$, 
one could not get a contrast-independent
condition number for the preconditioned system.
On the other hand,
when multiscale basis functions (one basis function per node) is 
used, the number of iterations and the condition number 
of the preconditioned system increase as the contrast increases. 
These results indicate that the preconditioner is optimal when the
spectral basis functions are used and the coarse spaces
include eigenvectors corresponding to important eigenvalues.
 Moreover, we observe that the dimension of the coarse space is smaller
when multiscale finite element basis functions are used an 
initial partition of unity.

\begin{table}[htb]
\centering
\footnotesize
\begin{tabular}{|l|r|r|r|r|r|r|}\hline
$\eta$     &         MS                                   &  GMsFEM with Lin $N_{rb}=3 $                & GMsFEM with MS $N_{rb}=3$          \\\hline
 $10^4$ & $125(6.52e+2)$       &     $ 45(22.79) $  &     $ 42 (15.16)$     \\
 $10^6$ & $260(6.14e+4)$       &     $ 37(8.94) $  &     $ 44 (13.15)$     \\
\hline\hline
 Dim & 81 (0.8\% of fine DOF)  & 862 (8.4\% of fine DOF) &  744 (7.4\% of fine DOF)   \\
\hline\hline 
\end{tabular}
\caption{Number of iterations until convergence 
and estimated condition number for  the PCG  and different values of 
the contrast $\eta$ with $\mu=0.5$. We set 
the tolerance to 1e-10. Here $H=0.1$ 
with $h=0.01$.}
\label{tab:newperm2-0ANI}
\end{table}

In Table \ref{tab:newperm1ms-10to6_1}, we present numerical results
to study the errors of GMsFEM when the contrast is $\eta=10^6$.
As there are three distinct spatial fields in the space of conductivities,
we choose at least three realizations.
In Table~\ref{tab:newperm1ms-10to6_1}, 
we compare the errors obtained by GMsFEM
when the online problem is solved with a corresponding
number of basis functions. 
We observe convergence with respect 
to the number of local eigenvectors.
The convergence with respect to  $N_{rb}$ 
can also be observed.  



\begin{table}
\centering\small 
\begin{tabular}{|l|r|r|r|r|||||||r|r|r|r|r|}\hline
$H=0.1$      & $N_{rb}=2$ & $N_{rb}=4$& $N_{rb}=4$\\\hline
 GMsFEM+0   &    $3.8(1204)$  &  $3.4(1360)$ & $3.2 (1517)$ \\
GMsFEM+1    &    $2.6(1325)$  &  $2.4(1481)$& $2.2(1638)$ )\\
GMsFEM+2   &  $2.1(1446)$  &  $1.8(1602)$& $1.7(1759)$\\\hline\hline 
\end{tabular}
\caption{Convergence results (energy norm in $\%$ and space dimension) with the increasing dimension
of the coarse space. 
Linear basis functions are used to generate the mass matrix.
Here, 
$h=0.01$, $\eta=10^8$, $\mu = 1/2$ (error with MsFEM $ 88\%$).}
\label{tab:newperm1ms-10to6_1}
\end{table}


\subsection{Some generalizations}
\label{sec:gen}

The  procedure proposed above can be applied to general linear problems
such as parabolic and wave equations.
The success of the method will depend on the local model reduction
that is encoded in $a^{\text{off}}_{\omega_i}$,  $a^{\text{on}}_{\omega_i}$,  
$s^{\text{off}}_{\omega_i}$,
and  $s^{\text{on}}_{\omega_i}$. These bilinear forms need to be 
appropriately defined for a given $L_\mu(\cdot)$.

One can apply GMsFEM to a linearized 
nonlinear problem where the operator is frozen at the current value
of the solution. In this case, one can treat the frozen value of
the solution as a scalar parameter on a coarse-grid block level. 
Note that if global model reduction techniques are used,
then one will have to  deal with a very large parameter space and these
computations will be prohibitively expensive.

To demonstrate
this concept, we assume that nonlinear equation is linearized
\[
L_\mu(u^{n+1};u^n)=f.
\]
For example, for the steady-state nonlinear quasilinear equation,
 we can consider the following linearization (see Section ~\ref{sec:ill_example})
\begin{equation}
\label{eq:linear}
-\mbox{div}(\kappa(x,u^n)\nabla u^{n+1})=f.
\end{equation}

To apply GMsFEM, one can consider
$u^n$ as a constant parameter, $\mu=\overline{u^n}$,
within each coarse-grid block. 
In this case, $\mu$ can be regarded as a parameter that represents
 the average of the solution in
each coarse block (see Section ~\ref{sec:ill_example}).
In the example of the steady-state Richards' 
equation, $u$ can be assumed to be a constant within a coarse-grid 
block. 
For general problems (e.g., if the linearization is around
$\nabla u$), one needs to use higher dimensional parameter space
to represent $\nabla u$ at a coarse-grid level. 

The construction of the offline space will follow the GMsFEM procedure
(see Section ~\ref{sec:ill_example}). 
We can construct a snapshot space either by taking local fine-grid functions
or
by using the value of the parameter corresponding to $\overline{u^j}$
that will appear in the linearization of the global system.
Furthermore, the offline space is constructed via a spectral
decomposition of the snapshot space as described
in  Section ~\ref{sec:ill_example}.

At the online stage, we consider an approximation of
 (\ref{eq:linear}) with $u^n$ replaced
by its average in each coarse-grid block. We denote this
approximate solution by $\widetilde{u}^{n+1}$
\[
-\mbox{div}(\kappa(x,\langle{\widetilde{u}^n}\rangle)\widetilde{u}^{n+1})=f,
\]
where $\langle{\widetilde{u}^n}\rangle$ 
is the average of $\widetilde{u}^n$ over the coarse regions.
For each parameter value $\mu$, which is
the average of the solution on a coarse-grid block
$\mu={1\over \omega_i} \int_{\omega_i}u$,  the
multiscale basis functions
are computed based on the solution of local problem. In particular,
for each $\omega_i$ and for each input parameter and 
$\mu={1\over \omega_i} \int_{\omega_i}u$,
 we formulate a quotient to find
 a subspace of $V_{\text{on}}^{\omega_i}$
where the space will be constructed for each $\mu$.
For the construction of the online space, we follow Section ~\ref{sec:ill_example}. The online space is computed by solving
an eigenvalue problem in $\omega_i$ using $V_{\text{off}}^{\omega_i}$
for the current value of $u^n$
see (\ref{eq:eigon}). Using dominant eigenvectors, we form
a coarse space as in (\ref{eq:G}) and solve the global coupled system
following (\ref{eq:globalG}) at the current $u^n$. For the numerical
example, we will consider the coefficient that has an affine
representation (see (\ref{hyp}))
which reduces the computational cost associated with
calculating the stiffness matrix in the online stage
(see page Section \ref{sec:RB}).

\begin{remark}[{Adaptivity in the parameter space}]

We note that one can 
use adaptivity in the parameter space to avoid computing
the offline space for a large range of parameters and compute the
offline space only for a short range of parameters and update the space.
This is,
in particular, the case for applications where one has
{\it a priori} knowledge about how the parameter enters into the problem.
To demonstrate this concept, we assume that the parameter space
$\Lambda$ can partitioned into a number of smaller parameter spaces
$\Lambda_i$, $\Lambda=\bigcup_i \Lambda_i$, where
$\Lambda_i$ may overlap with each other. Furthermore, the offline
spaces are constructed for each $\Lambda_i$. In the online stage,
depending on the online value of the parameter, we can decide
which offline space to use. This reduces the computational cost 
in the online stage. In many applications, e.g., in nonlinear problems,
one may remain in one of $\Lambda_i$'s for many iterations and
thus use the same offline space to construct the online space.

\end{remark}

We present a numerical example for
\begin{eqnarray*}
-\mbox{div}(\lambda(x,u)\nabla u ) = f, \label{eq:variational form2}
\end{eqnarray*}
with $u=0$ on $\partial D$ and
\begin{equation}
\lambda(x,u) = \lambda_0(x,u)\left( \kappa_1(x)+e^{\alpha u}\kappa_2(x)\right),
\end{equation}
where $\kappa_1(x)$ and $\kappa_2(x)$ are defined as in the previous example
(see Figure \ref{illus-onedim}).
The main objective of this example is to demonstrate that one can
use $u$ in $\lambda(x,u)$  as a scalar parameter
within a coarse-grid block. 
In contrast, in global methods, 
if $u$ in $\lambda(x,u)$ is used as a parameter, 
it will be a high dimensional parameter.
We will take the average of 
$u$ as a parameter in each coarse-grid block
as discussed above.

In Table \ref{tab:newperm1ms-10to6RBNL}, we present numerical results 
to study the accuracy of GMsFEM. We take $f=1$.
In our numerical simulations, we use $10$ values for the averaged
solution in each block and solve local eigenvalue problems for 
$8$ dominant eigenvectors to construct
the snapshot space. From this space $(80$ snapshots in each $\omega_i$), 
we construct the offline space by selecting the dominant eigenvectors.
In Table \ref{tab:newperm1ms-10to6RBNL}, we present numerical
results and show the dimension of the online space at the last iteration.
From our numerical results, we observe that the errors are small
and decrease as we increase the dimension of the local spectral spaces.
The number of iterations needed to converge is the same ($3$) for all
coarse space dimensions.
On the other hand, the error is large if MsFEM with one basis function
per node is used. This error in the energy norm is $47\%$.

\begin{table}
\centering\small 
\begin{tabular}{|l|r|r|r|r|r|}\hline
coarse dim &  $\Lambda_*$ &$|\cdot|_{L^2}^2$& $|\cdot|_A^2$\\\hline
293   & 78.7 & $2.87$ &  $14.17$  \\
352  &  244 &$0.4$ &  $7.02$ \\
620   & 981 &$0.05$ & $2.9$ \\\hline\hline
\end{tabular}
\caption{Relative errors in energy norm  and the coarse space
dimension in the last iteration. 
Here, 
$\eta=10^4$ .}
\label{tab:newperm1ms-10to6RBNL}
\end{table}


\section{Conclusions}

In this paper, 
we propose a multiscale framework, the Generalized Multiscale
Finite Element Method (GMsFEM), for solving PDEs with multiple scales. 
The main objective is to propose a 
framework that extends MsFEMs to more general problems with complex 
input space that includes parameters, high contrast, and right-hand-sides 
or boundary conditions. 
The GMsFEM starts with a family of snapshots for the local
solutions. These snapshots can usually be  generated based on the solutions
of local problems or simply taking local fine-grid functions.
 First, based on the local snapshot space,
the offline space is constructed. The construction of the offline space
involves solving a spectral problem in the snapshot space.
This process introduces a prioritization on the snapshot space
across the input space. In the online stage of the simulations,
for each new parameter and a source term, the online multiscale basis
functions are constructed efficiently. We discuss various constructions.
For example, in the absence of parameters, there is no computational
work needed in the online stage. When the solution nonlinearly
depends on the input space parameters, 
the construction of the online coarse spaces involves solving a spectral
problem over the offline space. We also discuss the online correction
of the reduced solution via two-level domain decomposition methods.
The optimality of the preconditioners is demonstrated through a few
examples.
We show that the GMsFEM covers some of existing multiscale methods. 
The generalization of the GMsFEM to
nonlinear problems is also considered. 
We illustrate these methods through a few numerical examples.
Numerical examples  suggest that the proposed framework 
can be effective in studying multiscale problems with an input 
space dimension and multiple right-hand-sides.

\section*{Acknowledgements}

We would like to thank Ms. Guanglian Li for helping us with the computations
and providing some computational results.
Y. Efendiev's work is
partially supported by the
DOE and NSF (DMS 0934837 and DMS 0811180).
J.Galvis would like to acknowledge partial support from DOE.
 This publication is based in part on work supported by Award
No. KUS-C1-016-04, made by King Abdullah University of Science
and Technology (KAUST).

\bibliographystyle{plain}

\def\cprime{$'$}

\end{document}